\documentclass[11pt]{amsart}
\usepackage[mathscr]{eucal}
\usepackage{mathrsfs}
\usepackage{amsfonts}
\usepackage{amsmath}
\usepackage{amsthm}
\usepackage{amssymb}
\usepackage{latexsym}
\usepackage[all]{xy}
\usepackage[latin1]{inputenc}
\usepackage{amscd, amssymb}
\usepackage{upgreek}
\usepackage{pstricks, pst-node}
\usepackage{bm}
\theoremstyle{plain}
\newtheorem{defi}{Definition}[section]
\newtheorem{thm}[defi]{Theorem}
\newtheorem{lem}[defi]{Lemma}
\newtheorem{prop}[defi]{Proposition}
\newtheorem{coro}[defi]{Corollary}
\theoremstyle{definition}
\newtheorem{rem}[defi]{Remark}
\newcommand{\pf}{\textbf{Proof.}}
\newcommand{\epf}{\begin{flushright} $\Box$\end{flushright}}
\newcommand{\End}{\mathrm{End}}

\newcommand{\Hom}{\mathrm{Hom}}

\newcommand{\tr}{\mathrm{Tr}}

\newcommand{\dete}{\mathrm{det}}
\newcommand{\di}{\mathrm{dim}\,}

\newcommand{\Rep}{\mathrm{Rep}}

\newcommand{\U}{\mathrm{U}}

\newcommand{\C}{\mathbb{C}}

\newcommand{\N}{\mathbb{N}}

\newcommand{\Z}{\mathbb{Z}}

\newcommand{\B}{\mathcal{B}}
\newcommand{\E}{\mathcal{E}}
\newcommand{\ii}{\underline{i}}

\newcommand{\M}{\scriptscriptstyle{M}}

\newcommand{\chik}{\chi^{(k)}}
\newcommand{\chikti}{\tilde{\chi}^{(k)}}

\newcommand{\e}{\epsilon}
\newcommand{\G}{\Gamma}
\newcommand{\GG}{\mathbf{\G}_n}
\newcommand{\g}{\gamma}
\newcommand{\Gg}{\mathrm{G}}
\newcommand{\GL}{\mathrm{GL}}
\newcommand{\gln}{\mathfrak{gl}_N}
\newcommand{\SLt}{\mathrm{SL}_2}
\newcommand{\HH}{\mathcal{H}\hspace{-4.7pt}\mathcal{H}}
\begin{document}

%

\title{Representations of Gan-Ginzburg algebras}
\author{Silvia Montarani}
\email{silvia.montarani@utoronto.ca}
\address{Department Mathematical and Computational Sciences,
University of Toronto Mississauga, 3359 Mississauga Road North,
Mississauga, ON L5L 1C6}
\begin{abstract}
Given a quiver, a fixed dimension vector, and a positive integer $n$, we construct a functor from the category of $D$-modules on  the space of  representations of the quiver to the category of modules over a corresponding Gan-Ginzburg algebra of rank $n$. When the quiver is affine Dynkin we obtain an explicit construction of representations of the corresponding wreath-product symplectic reflection algebra of rank $n$. When the quiver is star-shaped, but not finite Dynkin, we use this functor to obtain a Lie theoretic construction of representations of a ``spherical" subalgebra of the Gan-Ginzburg algebra  isomorphic to a rational generalized double affine Hecke algebra of  rank $n$. Our functors are a generalization of the type $A$ and type $BC$ functors from \cite{CEE} and \cite{EFM} respectively.   
\end{abstract}
\maketitle
\section{Introduction}

The present paper is a step in a program, proposed by Etingof, Loktev, Oblomkov, and Rybnikov,  that aims to connect the representation theory of Gan-Ginzburg algebras and symplectic reflection algebras with Lie theory  (see \cite{ELOR}).  

In \cite{CEE} Calaque, Etingof, and Enriquez give a Lie theoretic construction of some  representations of the  Cherednik algebra of type $A_{n-1}$ (see \cite{CEE}, Section 9.6). Specifically,  for any pair of positive integers $N$ and $n$, and any $D$-module $M$ over $\mathfrak{sl}_N$, they define an action of the Cherednik algebra of type $A_{n-1}$ with parameter $N/n$ on the space  
$$F_n(M):=(M\otimes (\C^N)^{\otimes n})^{\mathfrak{sl}_N},$$  
where the action of $\mathfrak{sl}_N$ on $M$ is by adjoint vector fields. 

In the same paper,  the authors extend their results  to  the case of the trigonometric degenerate double affine Hecke algebra (dDAHA) of type $A_{n-1}$, of which the Cherednik algebra is a further degeneration (see \cite{CEE}, Section 9.14).  For this purpose,  $D$-modules on the Lie algebra $\mathfrak{sl}_N$ are substituted with $D$-modules  on the group $\mathrm{SL}_N$. The construction for Cherednik algebras can be obtained from the one for dDAHA by rational degeneration.

In \cite{EFM} Etingof, Freund, and Ma  generalize the construction of \cite{CEE}  from  type $A$  to  type $BC$. Namely, the authors define a functor from the category of twisted $D$-modules on the symmetric space $\mathrm{GL}_N/\mathrm{GL}_p\times\mathrm{GL}_q$ (where $p+q=N$) to the category of representation of a corresponding dDAHA of type $BC_n$ (see \cite{EFM}, Section 5).

The  main result of this paper is a generalization of the construction of \cite{CEE} for Cherednik algebras to the  case of   arbitrary symplectic reflection algebras of wreath product type of any given rank $n$. These are symplectic reflection algebras attached to a  reflection group of the form $ \Gamma^n\rtimes S_n$, where $\Gamma$ is a finite subgroup of $\mathrm{SL}_2$, and to its natural representation on the space $\C^{2n}$. In particular, when $\Gamma$ is the trivial group, the corresponding wreath product symplectic reflection algebra of rank $n$ is the Cherednik algebra of type $A_{n-1}$, while when $\Gamma$ is the cyclic group of order two, $\Z/2\Z$,  it is the Cherednik algebra of type $BC_n$.  The generalization  of the functor $F_n$ that we propose was first suggested by Victor Ginzburg (see \cite{CEE}, Section 9.16) and it passes through the  representation theory of Gan-Ginzburg algebras, which can be seen as the link between the representation theory of symplectic reflection algebras and the theory of quiver representations.

Gan-Ginzburg  algebras are  one-parameter deformations of the wreath product of a symmetric group with the deformed preprojective algebra of a quiver. 
The  deformed preprojective algebra $\Pi_{\lambda}(Q)$ of a quiver $Q$, introduced by Crawley-Boevey and Holland in \cite{CBH}, is a deformation of the path algebra of the double quiver $\overline{Q}$ (the quiver obtained  by adding for each arrow in $Q$ an arrow in the opposite direction) with parameter $\lambda\in \C^I$, where $I$ is the set of vertices of $Q$. The representation theory of deformed preprojective algebras was studied in \cite{CBH}. The rank $n$ Gan-Ginzburg algebra $\mathcal{A}_{n,\lambda, \nu}(Q)$ is a deformation of the semi-direct product $\Pi_{\lambda}(Q)^{\otimes n}\rtimes\C[S_n]$,  depending
 on a single complex parameter $\nu$. In many important cases this deformation has the PBW (flatness) property.  
 
The original motivation  for the introduction of Gan-Ginzburg algebras comes from the representation theory of symplectic reflection algebras of wreath product type. Specifically  (see \cite{GG}) when $Q$ is an affine Dynkin quiver the algebra  $\mathcal{A}_{n,\lambda, \nu}(Q)$ is Morita equivalent to a  symplectic reflection algebra attached to the wreath product $\Gamma^n\rtimes S_n$, where  $\Gamma$  is associated to $Q$  via the classical  McKay correspondence (in the rank $1$ case this Morita equivalence was proved in \cite{CBH}).  

In this paper, for a fixed dimension vector $\alpha=\left\{\alpha_i\right\}_{i\in I}$ of the quiver $Q$, we consider the representation space $\mathrm{Rep}_{\alpha}(Q)$. The cotangent bundle of this  space is identified with the symplectic vector space $\mathrm{Rep}_{\alpha}(\overline{Q})$, and it is endowed with a symplectic action of the Lie group of basis changes $\Gg(\alpha)=\prod_{i\in I}\mathrm{GL}_{\alpha_i}$, and the corresponding infinitesimal action of its Lie algebra $\mathfrak{g}(\alpha)=\prod_{i\in I}\mathfrak{gl}_{\alpha_i}$. This gives a map (Weil representation) from the Lie algebra $\mathfrak{g}(\alpha)$ to the   algebra $W_{\alpha}=\mathcal{D}(\mathrm{Rep}_{\alpha}(Q))$ of differential operators with polynomial coefficients on the space   $\mathrm{Rep}_{\alpha}(Q)$. For a fixed $n$, consider the vector space $U^{\otimes n}$, where $U=\oplus_{i\in I}\C^{\alpha_i}$. This space carries natural actions of the Lie algebra $\mathfrak{g}(\alpha)$ and of the group $S_n$. Let $M$ be any $W_{\alpha}$-module,  and  let  $\chi$ be a character of $\mathfrak{g}(\alpha)$ that agrees on the scalars with the action of $\mathfrak{g}(\alpha)$ on $U^{\otimes n}$. We define an $\mathcal{A}_{n,\lambda,-1}$-module structure  on the $\chi$-equivariant space 
$$
F_{n,\chi}(M)=(M\otimes U^{\otimes n})^{\mathfrak{g}(\alpha)}_{\chi},
$$ 
that is the vector subspace of $M\otimes U^{\otimes n}$ on which $\mathfrak{g}(\alpha)$ acts by the character $\chi$. In this construction the parameter $\lambda$ is a certain explicit function of $\alpha$ and $\chi$.  We thus get a functor $F_{n,\chi}$ from the category of $W_{\alpha}$-modules to the category of representations of $\mathcal{A}_{n,\lambda,-1}(Q)$. When $Q$ is the Jordan quiver $\tilde{A}_0$ (the quiver with one vertex and one loop) our construction reduces to the one of \cite{CEE}, while when  $Q$ is the cyclic quiver of length two (that is the quiver with two vertices and two arrows pointing in opposite directions) it reduces to the rational limit of the functor from \cite{EFM}.

Another interesting aspect of Gan-Ginzburg algebras is their relation with rational generalized double affine Hecke algebras (rational GDAHA).  Rank one GDAHA attached to star-shaped  affine Dynkin diagrams were introduced by Etingof, Oblomkov and Rains in \cite{EOR} (except in the $D_4$ case for which they were already known thanks to the  work of Sahi and Stockman \cite{Sa,St}).  In \cite{ER}  Etingof and Rains extended the definition of GDAHA   to any  star-shaped graph that is not a finite Dynkin diagram. Rank one GDAHA are deformations of group algebras of appropriate Coxeter groups and have the PBW property. Later, in \cite{EGO},  Etingof, Gan, and Oblomkov  introduced GDAHA of higher rank attached to any star-shaped graph that is not a finite Dynkin diagram. GDAHA of rank $n$  are deformations of the wreath product of a rank one GDAHA with the symmetric group $S_n$, and are quotients of group algebras of appropriate braid groups. It is shown in \cite{EOR} that rank one GDAHA attached to affine Dynkin star shaped diagrams provide quantizations of del Pezzo surfaces (with a singular genus one curve removed).   Higher rank GDAHA are expected to give quantizations of spaces of Calogero-Moser type, which are some deformations of Hilbert Schemes of such surfaces. In \cite{EGO} the authors introduced a degenerate version of GDAHA called rational GDAHA (see \cite{EGO}, Section 2.2). In the affine case, and for values of the deformation parameters lying on a special hyperplane but otherwise generic, the rational GDAHA is finite over its center that can be identified with the algebra of regular functions on a corresponding Calogero-Moser space. 

 In \cite{EGO}  the authors prove that, when $Q$ is a star-shaped quiver that is not finite Dynkin,  the algebra $\mathcal{A}_{n,\lambda,\nu}(Q)$ contains a  subalgebra isomorphic to a corresponding rational GDAHA. This is the ``spherical" subalgebra $e_0^{\otimes n}\mathcal{A}_{n,\lambda,\nu}(Q)e_0^{\otimes n}$, where $e_0$ is the idempotent in the path algebra of $\overline{Q}$ that corresponds to the nodal vertex of the star-shaped quiver $Q$ (the nodal vertex is the vertex from which the legs of the star originate).  Thus, in particular, by composing our functor $F_{n, \chi}$ with the projector $e_0^{\otimes n}$ we get a  functor $e_0^{\otimes n}F_{n,\chi}$ from the category of $W_{\alpha}$-modules to the category of modules over the rational GDAHA. 

In the second part of our paper we  give an explicit Lie theoretic construction of some finite dimensional representations of  rational GDAHA of higher rank.  Specifically, suppose $D$ is a star-shaped graph (but not finite Dynkin) with $m$ legs whose lengths (including the node) are $\ell_1,\dots,\ell_m$ respectively.  Fix a  positive integer $N$ and  a character $\chi'$ of $\mathfrak{gl}_N$. For every index $i=1,\dots,m$, let $V_i$ be a $\ell_i$-stepped irreducible finite dimensional $\mathfrak{gl}_N$-module, that is an irreducible representation  such that the tensor product $V_i\otimes\C^N$ with the vector representation has at most $\ell_i$ distinct irreducible components.  For every   $n>1$, we define a representation of the rational GDAHA of rank $n$ attached to the graph $D$ on the $\chi'$-equivariant subspace 
$$
E_{n, \chi'}(V_1,\dots, V_m):=(V_1\otimes\dots\otimes V_m\otimes (\C^N)^{\otimes n})^{\mathfrak{gl}_N}_{\chi'}.
$$
We prove that  we can obtain some of these representations by restriction from the representations we constructed for Gan-Ginzburg algebras. In other words, we show that some of these representations are in the image of the functor  $e_0^{\otimes n}F_{n,\chi}$.    
Our main tool is a result from   \cite{ELOR}  establishing an isomorphism between the  algebra of twisted differential operators on a space of representations  for the quiver of type $A$ (where all arrows are oriented in the same direction and  the coordinates of the dimension vector are nonzero and strictly increasing in the direction of the arrows) and the algebra of twisted differential operators on a corresponding partial flag variety.
Namely, let $D$ be  a star shaped diagram that is not finite Dynkin, and let $Q$ be the quiver obtained from $D$ by assigning to all edges the orientation toward the node. Let $\alpha$ be a dimension vector such that $\alpha_0=N$ and the other components  are all nonzero and strictly increasing  along each leg moving toward the node.  We apply the functor $F_{n,\chi}$ to the $W_{\alpha}$-module $M=\C[\mathrm{Rep}_{\alpha}(Q)]$  of polynomial functions on $\mathrm{Rep}_{\alpha}(Q)$.  Projecting via the idempotent $e_0^{\otimes n}$  we get a representation of the corresponding  rational GDAHA on the space  
$$e_0^{\otimes n}F_{n,\chi}(M)=(M\otimes e_0^{\otimes n}U^{\otimes n})^{\mathfrak{g}(\alpha)}_{\chi}=(M\otimes(\C^N)^{\otimes n})^{\mathfrak{g}(\alpha)}_{\chi}.$$ 
This equivariant space can be obtained in two steps:  we can  first take the equivariant space with respect to the action of the Lie algebras of the groups of basis changes at the non-nodal vertices, and then the equivariant subspace with respect to the Lie algebra $\mathfrak{gl}_N$ of the group $\mathrm{GL}_N$ of basis changes at the nodal vertex. By the above mentioned result from \cite{ELOR} applied to the legs of $Q$ (that are quivers of type $A$), and for an appropriate choice of character $\chi$, after the first step we get a space of the form  $V_1\otimes\dots\otimes V_m\otimes (\C^N)^{\otimes n}$, where each vector space   $V_i$ is a finite dimensional irreducible $\ell_i$-stepped $\mathrm{gl}_N$-module whose highest weight is an explicit function of  $\chi$ and $\alpha$.
 
  We expect  our construction for rational GDAHA to have a generalization to the case of non-degenerate GDAHA, involving quantum $D$-modules and quantum groups. Such a generalization was already obtained by Jordan in the type $A$ case (see \cite{J}), and by Jordan and Ma  in the type  $BC$ case (see \cite{JM}). 

The paper is organized as follows. In Section \ref{ggaldefi} we recall the definition and main properties of Gan-Ginzburg algebras, as well as their relation to wreath-product symplectic reflection algebras. In Section \ref{ggrep} we construct the functor $F_{n,\chi}$ and compare it with the functors from \cite{CEE} and \cite{EFM} in the special case of the Cherednik algebra of type $A_{n-1}$ and $BC_{n}$ respectively. In Section \ref{gdahadefisec} we recall the definition of the rational GDAHA and its relation to Gan-Ginzburg algebras. In Section \ref{gdahasrep} we give a Lie theoretic construction of representations for rational GDAHA. In Section \ref{typeAqhr} we apply the results of \cite{ELOR} to obtain a Borel-Weil type construction of some finite dimensional irreducible $\mathrm{gl}_N$-modules starting from  modules over  algebras of differential operators on   spaces of representations of  quivers of type $A$. In Section \ref{comparison} we use the results of Section \ref{typeAqhr} to compare the representations of Gan-Ginzburg algebras and rational GDAHA we constructed in Section \ref{ggrep} and in Section \ref{gdahasrep} respectively.

 \subsection{Acknowledgments.} I  am very grateful to Pavel Etingof, Sergey Loktev,  Alexey Oblomkov, and Leonid Rybnikov for involving me in their project and for useful discussions. This paper was written while I was a Postdoctoral Fellow at the University of Toronto. I would like to thank my sponsors Sergey Arkhipov, James Arthur, Valentin Blomer, and  Eckhard  Meinrenken for their support.

\section{Gan-Ginzburg algebras}\label{ggaldefi}
\subsection{Definition of the Gan-Ginzburg algebra.} We recall some definitions from \cite{GG}. As in the previous section, let $Q$ be a connected quiver, and let $I$ be the set of vertices of $Q$. The \textit{double} quiver $\overline{Q}$ of $Q$ is the quiver obtained from $Q$ by adding for each arrow $i\stackrel{a}{\longrightarrow} j$ in $Q$ a reverse arrow $j\stackrel{a^{\ast}}{\longrightarrow}i$. 
In general, for an arrow $i\stackrel{a}{\longrightarrow} j$ we denote by $t(a)=i$ its tail (or source), and by $h(a)=j$ its head (or target).

Let $B:=\bigoplus_{i\in I}\C\,e_i$ be the finite dimensional semisimple algebra over $\C$ with basis formed by orthogonal idempotents $\left\{e_i\right\}_{i\in I}$ corresponding to the vertices of $\overline{Q}$ (or equivalently  of $Q$). Let $E$ be the  $\C$-vector space with basis formed by the set of arrows $\left\{a\in \overline{Q}\right\}$. Then $E$ is naturally a $B$-bimodule, and as such it decomposes as $E=\bigoplus_{i,j\in I} E_{ij}$, where $E_{ij}=e_iEe_j$ is the submodule spanned by the arrows $a\in \overline{Q}$ with $h(a)=i$ and $t(a)=j$. The path algebra of $\overline{Q}$ is $\C\overline{Q}:=T_BE=\bigoplus_{n\geq 0}T^n_BE$, where $T^n_BE=E\otimes_B\cdots\otimes_BE$ is the $n$-fold tensor product.  The algebra  $\C\overline{Q}$ can be equivalently described as the algebra of paths on the oriented graph $\overline{Q}$, where each idempotent $e_i$ is identified with the trivial path starting and ending at the corresponding vertex $i$, and the composition of paths that are not compatible is zero.

Let  $n$ be  a positive integer and set $\B:=B^{\otimes n}$. Let $\ell\in[1,n]$ be an index. We define the $\B$-bimodules $\E_{\ell}$ and $\E$ as follows:
$$
\E_{\ell}:=B^{\otimes (\ell-1)}\otimes E\otimes B^{\otimes (n-\ell)}\qquad  \qquad \E:=\bigoplus_{1\leq \ell \leq n}\E_{\ell}.
$$
We can thus form the tensor algebra $T_{\B}\E$. If we denote by $\ii=(i_1,\dots, i_n)\in I^n$ a multi-index, the unity in $T_{\B}\E$ is given by the sum of idempotents $\sum_{\ii\in I} e_{i_1}\otimes e_{i_2}\otimes\cdots\otimes e_{i_n}$.
The symmetric group $S_n$ acts naturally on $\E$. This action extends to $T_{\B}\E$, and we can form the smash product algebra $T_{\B}\E\rtimes \C[S_n]$.

For any path $p\in \C\overline{Q}$ and $\ell\in [1,n]$, consider the element $(p)_{\ell}\in T_{\B}\E\rtimes\C[S_n]$ defined by the following formula
$$
(p)_{\ell}:=1\otimes\cdots\otimes 1\otimes p\otimes 1\otimes\cdots\otimes 1
$$
where $p$ is placed in the $\ell$-th position.

Let $\nu$ be a complex number, and let  $\lambda=\sum_{i\in I}\lambda_ie_i$ be an element of $B$.
For any pair of distinct indices $\ell, m\in [1,n]$, let $s_{\ell m}\in S_n$ be the transposition  $\ell\longleftrightarrow m$. Define the projector $P_{\ell,m}\in T_{\B}\E\rtimes\C[S_n]$ as
$$
P_{\ell,m}=\sum_{i\in I}(e_i)_{l}(e_i)_m\in \B.
$$

The following definition is equivalent to \cite{GG}, Definition 1.2.3

\begin{defi}\label{ggsimple}
The Gan-Ginzburg algebra $\mathcal{A}_{n,\lambda, \nu}(Q)$ is the quotient of $T_{\B}\E\rtimes\C[S_n]$ by  the following relations.
\begin{enumerate}
\item[(I)] For any $\ell\in[1,n]$ 
$$
\sum_{a\in Q}[(a)_{\ell}, (a^{\ast})_{\ell}]=(\lambda)_{\ell}+\nu\sum_{m\neq \ell} P_{\ell, m}s_{\ell,m}.
$$
\item[(II)] For any $a, b\in \overline{Q}$, and any $\ell, m\in [1,n]$ with $\ell\neq m$
$$
\left[(a)_{\ell},(b)_{m}\right]=\left\{\begin{array}{rl}\nu(e_{h(a)})_{\ell}(e_{t(a)})_ms_{\ell m} &\mbox{if\ }  a=b^{\ast}\mbox{\ and\ } b\in Q\\
-\nu(e_{h(a)})_{\ell}(e_{t(a)})_ms_{\ell m} &\mbox{if\ }  b=a^{\ast}\mbox{\ and\ } a\in Q\\
0\phantom{(e_{h(a)})_{\ell}(e_{t(a)})_ms_{\ell m} } &\mbox{else}\end{array}\right.
$$
\end{enumerate}
\end{defi}
 
 The proof of the equivalence of the two definitions  is trivial. Indeed  relations (I) and (II) in Definition \ref{ggsimple} can be obtained from relations (i) and (ii) in  Definition 1.2.3 of \cite{GG} respectively, by summing over the indices $i_1,\dots, i_n$. Conversely, the above mentioned relations  (i) and (ii) can be obtained  by multiplying (I) and (II) on both sides by the idempotents $e_{i_1}\otimes\cdots\otimes e_{i_n}$ and $e_{i_1}\otimes\cdots\otimes\underbrace{e_{h(a)}}_{\ell}\otimes\cdots\otimes\underbrace{e_{t(a)}}_m\otimes\cdots\otimes e_{i_n}$ respectively. 
Thus the two sets of relations generate the same ideal in $T_{\B}\E\rtimes\C[S_n]$.

We recall that the algebra $\mathcal{A}_{n,\lambda,\nu}(Q)$ does not depend on the orientation of $Q$. Moreover a simultaneous rescaling of the parameters by a nonzero complex number does not change the algebra $\mathcal{A}_{n,\lambda,\nu}(Q)$ up to isomorphism.
 
\subsection{Relation to deformed preprojective algebras and  PBW property.}
It is straightforward to check that:
\begin{itemize}
\item[(i)] The algebra $\mathcal{A}_{1,\lambda,\nu}(Q)=\mathcal{A}_{1,\lambda}(Q)$ is independent of $\nu$  and is the deformed preprojective algebra 
$\Pi_{\lambda}(Q)$ of \cite{CBH} which is, by definition, the quotient of the path algebra $\C\overline{Q}$ by the following relation
$$
\sum_{a\in Q}[a,a^{\ast}]-\sum_{i\in I}\lambda_ie_i=0.
$$

\item[(ii)] $\mathcal{A}_{n,\lambda,0}(Q)=\prod_{\lambda}(Q)^{\otimes n}\rtimes \C[S_n]$.
\end{itemize}\medskip

It follows that $\mathcal{A}_{n,\lambda,\nu}(Q)$ is a  one-parameter deformation of $\Pi_{\lambda}(Q)^{\otimes n}\rtimes \C[S_n]$ with parameter $\nu$.  This deformation is not always flat. However one has the following theorem due to Gan and Ginzburg (see \cite{GG}, Theorem 2.2.1 and Remark 2.2.6).

Consider the filtration on $\mathcal{A}_{n,\lambda,\nu}(Q)$ obtained by assigning degree zero to the elements of $S_n$ and of $\mathcal{B}$, and degree one to each arrow $a\in \overline{Q}$. It is easy to see that there is a surjective homomorphism $\phi:\mathcal{A}_{n,0,0}(Q)\rightarrow \mathrm{gr}(\mathcal{A}_{n,\lambda,\nu})$. 
\begin{thm}[\cite{GG}]
If $Q$ is a connected quiver that is not of finite Dynkin type and  has no edge-loop then $\phi$ is an isomorphism.
\end{thm}

\epf

In this case we say that $\mathcal{A}_{n,\lambda,\nu}(Q)$ is a PBW  deformation of $\mathcal{A}_{n,0,0}=\prod_0^{\otimes n}(Q)\rtimes \C[S_n]$.

\subsection{Symplectic reflection algebras for wreath products.} The original motivation to study Gan-Ginzburg algebras comes from the theory of symplectic reflection algebras introduced by Etingof and Ginzburg in \cite{EG}.  In fact, in the case when $Q$ is an extended Dynkin quiver, Gan-Ginzburg algebras have been a fundamental tool in the study of the representation theory of the wreath-product symplectic reflection algebras (see \cite{G}, \cite{EM}, \cite{M}), of which we are going to recall the definition below.

Let $L$
be a $2$-dimensional complex vector space with a symplectic form
$\omega_L$, and consider the space $V=L^{\oplus n}$,
endowed with the induced symplectic form
$\omega_V={\omega_L}^{\oplus n}$. Let $\Gamma$ be a finite
subgroup of $Sp(L)$. Let $S_n$ act
on $V$ by permuting the factors. The wreath product group
$\GG:=\Gamma^n\rtimes S_n \subset Sp(V)$ acts naturally on
$V$. In what follows we write $\gamma_i\in \G_n$  for any
element $\gamma\in \Gamma$ seen as an element in the $i$-th factor
$\Gamma$ of $\GG$. 
It is clear that choosing a symplectic basis $x$, $y$ for $L$ we
can identify $\Gamma$ with a subgroup of $\SLt$. We denote by $x_i$, $y_i$ the corresponding vectors in the $i$-th summand $L$ of $V$. Similarly, given any $v\in V$ and $i\in [1,n]$, we write $v_i$  for  $v$ placed in the $i$-th summand $L$ of $V$. Let $t, k\in \C$, and let $c$ be a complex valued class function on $\G\smallsetminus\left\{1\right\}$. Write $c_{\g}=c(\g)$ for the value that $c$ assumes on an element $\g\in \G$. Denote by $TV$ the tensor algebra of $V$. The following presentation of the wreath product symplectic reflection algebra is due to Gan and Ginzburg (see \cite{GG}, Lemma 3.1.1).
\begin{defi} \label{dewre}
The symplectic reflection algebra $H_{t,k,c}(\GG)$ is the quotient
of the smash product $TV\rtimes\C[\GG]$ by the following relations:
\begin{itemize}
\item[{\emph{(R1)}}] 
For any $i\in [1,n]$: $$[x_{i}, y_{i}]
=  t- \frac{k}{2} \sum_{j\neq i}\sum_{\g\in\Gamma}
s_{ij}\g_{i}\g_{j}^{-1} + \sum_{\g\in\Gamma\backslash \{1\}}
c_{\g}\g_{i}.$$
\item[{\emph{(R2)}}] 
For any $u,v\in L$ and $i\neq j$:
$$[u_{i},v_{j}]= \frac{k}{2} \sum_{\g\in\Gamma} \omega_{L}(\g u,v)
s_{ij}\g_{i}\g_{j}^{-1}.$$
\end{itemize}
\end{defi}

We recall that the algebra $H_{t,k,c}(\GG)$ does not change up to isomorphism if the parameters $t,k,c$ are simultaneously rescaled by a nonzero constant; in particular if $t\neq 0$   we have that $H_{t,k,c}(\GG)\simeq H_{1,k/t,c/t}(\GG)$. Note that our parameter $k$ differs from the one in \cite{GG} by a sign. As we will see in Section \ref{sectioncherednik}, we make this choice so that, in the special case of the rational Cherednik algebra, our notation is consistent with the notation of \cite{EG} and \cite{E}, Example 7.5.

\subsection{Morita equivalence.} We will now recall the classical McKay correspondence. For a finite subgroup $\G\subset \SLt$, we denote by $\left\{N_i\right\}_{i\in I}$ a complete set of pairwise non-isomorphic representations. We consider the quiver with vertex set $I$, and such that the number of arrows from  $i$ to $j$ is the multiplicity of $N_i$ in $L\otimes N_j$, where $L$ is the defining representation of $\G$. In other words we have that $\dim E_{ij}=\dim \mathrm{Hom}_{\G}(N_i, N_j\otimes L)$, and since $L$ is a self-dual representation it follows that $\dim E_{ij}=\dim E_{ji}$. The quiver so obtained is the double of an affine Dynkin quiver of type ADE. Forgetting the orientation, this construction establishes a bijection between the conjugacy classes of finite subgroups of $\SLt$  and the affine Dynkin diagrams of type ADE. The following result, due to Gan and Ginzburg (see \cite{GG}, Theorem 3.5.2), is the original motivation for the introduction of Gan-Ginzburg algebras.

\begin{thm}\label{ggmorita}
The algebra $H_{t,k,c}(\GG)$ is Morita equivalent to the  algebra $\mathcal{A}_{n,\lambda,\nu}(Q)$, where $Q$ is the McKay quiver for $\G$, $\lambda_i$ is the trace of $t\cdot 1+\sum_{\g\neq 1}c_{\g}\g$ on $N_i$, and $\nu=-\frac{k|\G|}{2}$.
\end{thm}

For $n=1$ the above result is due to Crawley-Boevey and Holland (see \cite{CBH}, Corollary 3.5).
\subsection{The special case of the rational Cherednik algebra.}\label{sectioncherednik}
 When $\G$ is trivial or a cyclic group the algebra in Definition \ref{dewre} is a rational Cherednik algebra, and the Morita equivalence of Theorem \ref{ggmorita} is actually an isomorphism.
 
In particular, when $\G$ is  trivial  there is no parameter $c$ and   the algebra in Definition \ref{dewre} is isomorphic to the Cherednik algebra of type $A_{n-1}$, that is the Cherednik algebra $H_{t,k}(S_n,\C^n)$  attached to the reflection group $S_n$  and to its natural permutation representation on $\C^n$ (where the parameters are the same as the ones in \cite{EG},  and in \cite{E}, Example 7.5). It is also easy to check that in this case  the vectors $x_{i}-x_{i+1}, y_{i}-y_{i+1}$ where $i=1,\dots,n-1$ and the elements of $S_n$ generate a subalgebra isomorphic to the rational Cherednik algebra   $H_{2t,k}(S_n,\C^{n-1})$ attached to $S_n$ and to its reflection representation. In this case the McKay quiver of $\G$ is the Jordan quiver $\tilde{A}_0$ (that is the quiver with only one vertex and one loop) and  we have an isomorphism 
$$\mathcal{A}_{n,\lambda,\nu}(\tilde{A}_0)\simeq H_{\lambda,-2\nu}(S_n,\C^n).$$ 
It follows  that  $\mathcal{A}_{n,\lambda,\nu}(\tilde{A}_0)$ contains a subalgebra isomorphic to the rational Cherednik algebra  
$$H_{2\lambda, -2\nu}(S_n, \C^{n-1})\simeq H_{1, -\nu/\lambda}(S_n,\C^{n-1}).$$

Similarly, it is easy to see that when $\G$ is a cyclic group $\Z/m\Z$ of order $m\geq 2$  the algebra of Definition \ref{dewre} is isomorphic to the rational Cherednik algebra  $H_{t, k, c}((\Z/m\Z)^n\rtimes S_n,\C^n)$ attached to the complex reflection group $G(m,1,n)=(\Z/m\Z)^n\rtimes S_n$ and to its reflection representation on $\C^n$. In particular when $m=2$ we get the rational Cherednik algebra of type $BC_n$.
In this case the McKay quiver of $\G$ is the cyclic quiver of length two (that is the quiver with two vertices and two arrows pointing in different directions).  In this case the parameter $c$ is just a complex number, the parameter  $\lambda$ is a pair of complex numbers $(\lambda_1,\lambda_2)$, and there is an isomorphism
 $$
 \mathcal{A}_{n,\lambda,\nu}(Q)\simeq H_{t,k,c}((\Z/2\Z)^n\rtimes S_n, \C^n)
 $$ 
 where it is easy to calculate that, according to Theorem \ref{ggmorita}, the parameters $(t,k,c)$ and $(\lambda,\nu)$ are related by the following formulas
\begin{equation}\label{parametersBCn}
t=\frac{\lambda_1+\lambda_2}{2},\qquad\qquad
c=\frac{\lambda_1-\lambda_2}{2},\qquad\qquad
k=-\nu
\end{equation} 
\begin{rem}\label{degeneration}
It is not hard to check that the algebra $H_{t, k, c}((\Z/m\Z)^n\rtimes S_n,\C^n)$ can be obtained by rational degeneration from the (trigonometric) degenerate  double affine Hecke algebra (dDAHA) $\HH(t,k_1,k_2,k_3)$ described in \cite{EFM}, Section 3.1.  Roughly speaking, the rational degeneration is obtained by setting $X_i=e^{\hbar \tilde{x}_i}$ and $\tilde{y_i}=\hbar y_i$ where $X_i$ and $y_i$ are the generators of $\HH(t,k_1,k_2,k_3)$ in \cite{EFM},  and by taking the limit   $\hbar \rightarrow 0$. The  generators $\tilde{x_i}$, $\tilde{y_i}$ of the rational degeneration are related to the generators $x_i$, $y_i$ in Definition \ref{dewre} by the formulas $x_i=\tilde{y}_i$ and $y_i=\tilde{x}_i$. The relations between the parameters $(k,c)$ and $(k_1,k_2,k_3)$ are as  follows
\begin{equation}\label{parametersdegeneration}
k=2k_1,\qquad\qquad c=-(2k_2+k_3).
\end{equation}

\end{rem}

\section{ Construction of the functors $F_{n,\chi}$}\label{ggrep}

\subsection{The main theorem.} Let  $Q$ be a connected quiver, and  let $\alpha=\left\{\alpha_i\right\}_{i\in I}$ be a dimension vector for $Q$. We denote by $\Rep_{\alpha}(\overline{Q})$ the space of matrix representations of dimension vector $\alpha$ of the double quiver $\overline{Q}$ . We have
$$
\Rep_{\alpha}(\overline{Q})=\bigoplus_{i,j}\Hom_{\C}\left(E_{ij},\Hom_{\C}\left(\C^{\alpha_j},\C^{\alpha_i}\right)\right).
$$ 
We  denote a representation in $\Rep_{\alpha}(\overline{Q})$ by  $\rho=\left\{\rho(a)\right\}_{a\in \overline{Q}}$.

Observe that if we denote by $Q^{op}$  the opposite quiver, that is the quiver  obtained from $Q$ by reversing the orientation of the arrows, we have the following decomposition 
$$
\Rep_{\alpha}(\overline{Q})=\Rep_{\alpha}(Q)\oplus \Rep_{\alpha}(Q^{op})=\bigoplus_{a\in Q}\Hom_{\C}\left(\C^{t(a)}, \C^{h(a)}\right)\oplus\Hom_{\C}\left(\C^{h(a)}, \C^{t(a)}\right).
$$
Using the trace pairing, we obtain  canonical isomorphisms  
$$\Rep_{\alpha}(Q^{op})\simeq\Rep_{\alpha}(Q)^{\ast}$$ 
and   
$$
\Rep_{\alpha}(\overline{Q})\simeq \Rep_{\alpha}(Q)\times \Rep_{\alpha}(Q)^{\ast}\simeq T^{\ast}\Rep_{\alpha}(Q).
$$
It follows that $\Rep_{\alpha}(\overline{Q})$ is equipped with a natural symplectic form  that, using again the trace pairing, is identified with the form $\eta$ on $\Rep_{\alpha}(\overline{Q})$ defined by the following formula
$$
\eta(\rho,\tau)=\sum_{a\in Q}\tr\left(\rho(a)\tau(a^{\ast})-\tau(a)\rho(a^{\ast})\right)\qquad \forall\, \rho,\tau\in \Rep_{\alpha}(\overline{Q}).
$$
Let $W_{\alpha}$ be the Weyl algebra of the symplectic vector space $(\Rep_{\alpha}(\overline{Q}), \eta)$, that is the algebra generated by the vectors in $\Rep_{\alpha}(\overline{Q})$ with  defining relations
\begin{equation}\label{weylrel}
\left[\rho,\tau\right]=\eta(\rho, \tau) \qquad \forall\, \rho,\tau\in \Rep_{\alpha}(\overline{Q}).
\end{equation}
In other words $W_{\alpha}=\mathcal{D}(\Rep_{\alpha}(Q))$, where $\mathcal{D}(\Rep_{\alpha}(Q))$ denotes the algebra of differential operators with polynomial coefficients on the vector space $\Rep_{\alpha}(Q)$. Let $U=\oplus_{i\in I}\C^{\alpha_i}$ be the $B$-module in which, for any $i$, the idempotent $e_i$ acts as the projector on the subspace $\C^{\alpha_i}$.

For every $m\in \N$ we denote by $\GL_m$  the general linear group of  degree $m$ over $\C$, and by  $\mathfrak{gl}_m$   its Lie algebra. We set  $\Gg(\alpha)=\prod_{i\in I}\GL_{\alpha_i}$. Thus $\Gg(\alpha)$ is reductive, with Lie algebra  $\mathfrak{g}(\alpha)=\prod_{i\in I}\mathfrak{gl}_{\alpha_i}$, and it can be identified with the space of automorphisms of $U$ as a $B$-module. The natural action of $\Gg(\alpha)$  on the space $\Rep_{\alpha}(Q)$ by basis change automorphisms induces an action on the cotangent bundle $\Rep_{\alpha}(\overline{Q})$ that preserves the symplectic form $\eta$. It is a standard fact that, in this case, we have a natural homomorphism $\xi:\mathfrak{g}(\alpha)\longrightarrow W_{\alpha}$, known as the Weil representation. If for   $x\in \mathfrak{g}(\alpha)$ and $\rho\in \Rep_{\alpha}(\overline{Q})$ we denote by $x\rho$ the infinitesimal action, then the homomorphism $\xi$ is  given by the formula 
\begin{equation}\label{defixi}
\xi(x)=\frac{1}{2}\sum_p(x\rho^p)\rho_p
\end{equation}
where $\left\{\rho_p \right\}$, $\left\{\rho^p \right\}$ are dual bases of $\Rep_{\alpha}(\overline{Q})$ with respect to $\eta$, that is $\eta(\rho_p,\rho^q)=\delta_{pq}$.
Note that for  any $b\in \mathfrak{g}(\alpha)$ and $\rho\in \Rep_{\alpha}(\overline{Q})$ we have 
\begin{equation}
\left[\xi(x),\rho\right]=x\rho
\end{equation}
where $[\ ,\ ]$ denotes the commutator in the associative algebra $W_{\alpha}$.

It follows that every $W_{\alpha}$-module $M$ is  a $\mathfrak{g}(\alpha)$-module via the Weil representation. For every $n\in \N$ and every $W_{\alpha}$-module $M$ 
we can thus consider the $\mathfrak{g}(\alpha)$-module  $M\otimes U^{\otimes n}$ where the $\mathfrak{g}(\alpha)$-module structure on $U$  is given by the infinitesimal action induced by the natural action of $\Gg(\alpha)$.  

Let  $\chi$ be a character of $\mathfrak{g}(\alpha)$ that agrees on the scalars with the central character of the $\mathfrak{g}(\alpha)$-module $U^{\otimes n}$. We define the $\chi$-equivariant subspace  of $M\otimes U^{\otimes n}$ as
$$(M\otimes U^{\otimes n})^{\mathfrak{g}(\alpha)}_{\chi}:=\left\{w\in M\otimes U^{\otimes n}|xw=\chi(x)w,\,\;\forall\, x\in\mathfrak{g}(\alpha)\right\}.$$
Observe that, since scalars in $\Gg(\alpha)$ act trivially on $\Rep_{\alpha}(\overline{Q})$,  the map $\xi$ is zero on the scalars in $\mathfrak{g}(\alpha)$ and thus  scalars  act trivially on $M$. It follows that the  condition we imposed on  the character $\chi$ is a necessary condition for the corresponding equivariant space to be nonzero. 

For every choice of $\chi$ and $n$ as above we define a functor  $F_{n,\chi}:W_{\alpha}-\mathrm{mod}\longrightarrow \mathrm{Vect}$  by the formula
\begin{equation}\label{functordefi}
F_{n,\chi}(M)=(M\otimes U^{\otimes n})^{\mathfrak{g}(\alpha)}_{\chi} \qquad \forall\; M \mbox{\ \ in\ \ } W_{\alpha}-\mathrm{mod}.
\end{equation}

For each generator $(a)_{\ell}$ of the Gan-Ginzburg algebra $\mathcal{A}_{n,\lambda, \nu}(Q)$, where $\ell=1,\dots n,$ and $a\in \overline{Q}$,  we  define a linear operator on $F_{n,\chi}(M)$ that, for simplicity, we denote by the same symbol  
\begin{equation}\label{arrowmaps}
(a)_{\ell}:=\sum_{p}\rho_p|_{\M}(\rho^p(a))_{\ell}.
\end{equation}
Here  $\left\{\rho_p\right\}$ and $\left\{\rho^p \right\}$ are dual bases with respect to $\eta$, and the subscript $\ell$ in $(\rho^p(a))_{\ell}$
indicates that the operator $\rho^p(a)$ acts on the $\ell$-th factor of $U^{\otimes n}$. We also have an action of the algebra $\mathcal{B}\rtimes\C[S_n]=B^{\otimes n}\rtimes\C[S_n]$ on $F_{n,\chi}(M)$, induced by the natural action of $\mathcal{B}\rtimes\C[S_n]$ on $U^{\otimes n}$. Indeed, since $\mathfrak{g}(\alpha)=\End_{B}(U)$, the $\mathcal{B}$ action commutes with the $\mathfrak{g}(\alpha)$ action on $U^{\otimes n}$, and thus it preserves the $\chi$-equivariant  subspace  (\ref{functordefi}). Moreover
the $S_n$ action on $U^{\otimes n}$ commutes with the $\Gg(\alpha)$ action, thus with the $\mathfrak{g}(\alpha)$ action,  and so it also preserves the above mentioned subspace. The smash product relations of $\mathcal{B}\rtimes\C[S_n]$ are obviously satisfied. 

We recall now that the symmetrized Ringel form attached to $Q$ is the bilinear form on $\Z^{|I|}$ given by the formula
\begin{equation}\label{ringelsym}
(\alpha,\alpha')=2\sum_{i\in I}\alpha_i\alpha_i'-\sum_{a\in Q}\alpha_{t(a)}\alpha_{h(a)}'+\alpha_{h(a)}\alpha_{t(a)}' \qquad \forall\, \alpha,\alpha'\in \Z^{|I|}. 
\end{equation}
Let  $\left\{\epsilon _i\right\}_{i\in I}$ be the $\Z$-basis of $\Z^{|I|}$, where $\epsilon_i=(0,\dots,1_{j},\dots,0)$, and let $C=(c_{ij})$ be the matrix of the above form  in this basis. Then   $c_{ij}=2\delta_{ij}-a_{ij}$, where $A=(a_{ij})$ is the adjacency matrix of the quiver $\overline{Q}$, that is $a_{ij}=\di E_{ij}$ \footnote{ Note that here we adopt the definition of adjacency matrix for a directed graph in which the $i$-th diagonal entry is the number of loops  at the corresponding vertex $i$, and not twice such number. Thus, in our case, the $i$-th diagonal entry of the matrix $A$ is  the number of loops of $\overline{Q}$ at the vertex $i$, which is twice the number of loops of $Q$ at the same vertex.}. 

Observe that any character $\chi$ of $\mathfrak{g}(\alpha)$  can be written as a sum of the form
$$\chi=\sum_{i\in I}\chi_i\tr_{\mathfrak{gl}_{\alpha_i}}$$ where $\chi_i$ is in  $\C$. The condition that  $\chi$  agrees on the scalars with the central character of $U^{\otimes n}$ corresponds to the condition
\begin{equation}\label{conditionchis}
\sum_{i\in I}\chi_i \alpha_i-n=0
\end{equation}
on the complex numbers $\chi_i$ for $i\in I$. 

The following theorem is the first main result of this paper.
\begin{thm}\label{representationtheorem}
The operators $(a)_{\ell}$ and the action of $\mathcal{B}\rtimes\C[S_n]$ on $F_{n,\chi}(M)$ combine into an action of the Gan-Ginzburg algebra $\mathcal{A}_{n,\lambda,\nu}$ where
$$
\lambda_i=\chi_i-\frac{1}{2}\sum_{j}c_{ji}\alpha_j
$$
and
$$
\nu=-1.
$$
\end{thm}
\subsection{Proof of Theorem \ref{representationtheorem}.} In preparation for the proof of Theorem \ref{representationtheorem} we will now prove two auxiliary lemmas. The first lemma describes some properties of  dual basis for the symplectic form $\eta$ on $\Rep_{\alpha}(\overline{Q})$. These properties are straightforward but for the sake of clarity in our computations we prefer to present them here separately, and  use them freely in the rest of the paper.
\begin{lem}
Let $\left\{\rho_p\right\}$ and $\left\{\rho^p\right\}$ be dual basis with respect to $\eta$, and let $x$ be an element of $\mathfrak{g}(\alpha)$. Then
\begin{itemize}
\item[i)] $$\rho_q=\sum_p\eta(\rho_p,\rho_q)\rho^p=-\sum_p\eta(\rho_q,\rho_p)\rho^p;$$
\item[ii)] $$\sum_p (x\rho^p)\rho_p=\sum_p\rho_p(x\rho^p).$$
\end{itemize}
\end{lem}
\pf\
Part i) is trivial. Let us prove part  ii). We have
\begin{align*} 
 \sum_p (x\rho^p)\rho_p &=\sum_p\rho_p(x\rho^p)+\sum_p\eta(x\rho^p, \rho_p)\\ 
                        &=\sum_p\rho_p(x\rho^p)-\sum_p\eta(\rho^p, x\rho_p)\\
                        &=\sum_p\rho_p(x\rho^p)+\sum_p\eta(-\rho^p, x\rho_p)\\ 
                        &=\sum_p\rho_p(x\rho^p)+\sum_p(-\rho^p)(x\rho_p)-\sum_p(x\rho_p)(-\rho^p) 
\end{align*}
where we used the fact that $\eta$ is invariant with respect to the $\mathfrak{g}(\alpha)$-action. From the above calculations it follow that
\begin{equation}\label{xinvariance}
 \sum_p (x\rho^p)\rho_p+\sum_p(x\rho_p)(-\rho^p)=\sum_p\rho_p(x\rho^p)+\sum_p-\rho^p(x\rho_p).
\end{equation}
We observe now that the right hand side and the left hand side of the identity in ii) do not depend on the choice of the pair of dual basis, and that if  $\left\{\rho_p\right\}$ is dual to $\left\{\rho^p\right\}$  with respect to $\eta$ then $\left\{-\rho^p\right\}$ is dual to $\left\{\rho_p\right\}$. Thus formula \eqref{xinvariance} can be rewritten as 
$$
2\sum_p(x\rho^p)\rho_p=2\sum_p \rho_p(x\rho^p)
$$
which is the desired identity.

\epf

Let $a$, $b$ be arrows in $\overline{Q}$, and  consider  the following operator:
\begin{equation}\label{operatorrelationII}
F:=\sum_q\rho_q(a)\otimes\rho^q(b):U\otimes U\longrightarrow U\otimes U.
\end{equation}
Our second lemma describes the action of $F$ on $U\otimes U$. 
\begin{lem}\label{operatorflip} We have that
$$
F=\left\{\begin{array}{ll} -flip \circ \left(e_{t(a)}\otimes e_{h(a)}\right) &\mbox{if\ } a=b^{\ast}\mbox{\ and\ } b\in Q\\\phantom{-}flip \circ \left(e_{h(a)}\otimes e_{t(a)}\right) &\mbox{if\ } b=a^{\ast}\mbox{\ and\ } a\in Q \\\phantom{-}0 &\mbox{else\ }\end{array}\right..
$$
where by $flip$ we denote the linear map that exchanges the  factors in $U\otimes U$.
\end{lem}
\pf\ 
Observe  that if $E_{ij}$ denotes the vector space spanned by the arrows in the double quiver $\overline{Q}$ with tail $j$ and head $i$ then we can write $E_{ij}=E_{ij}(Q)\oplus E_{ij}(Q^{op})$, where  $E_{ij}(Q)$ and $E_{ij}(Q^{op})$ denote the vector spaces spanned by the arrows with tail $j$ and head $i$ in the quiver $Q$ and $Q^{op}$ respectively. Thus we have
\begin{align}\label{decompositionsymplectic}
\Rep_{\alpha}(\overline{Q})=&\bigoplus_{i,j}\Hom_{\C}\left(E_{ij},\Hom_{\C}\left(\C^{\alpha_j},\C^{\alpha_i}\right)\right)\\
       =&\bigoplus_{i,j}\Hom_{\C}\left(E_{ij}(Q),\Hom_{\C}\left(\C^{\alpha_j},\C^{\alpha_i}\right)\right)\oplus \Hom_{\C}\left(E_{ij}(Q^{op}),\Hom_{\C}\left(\C^{\alpha_j},\C^{\alpha_i}\right)\right).\nonumber
\end{align}    

For any pair $i,j\in I$ take $\left\{u^{ij}_{rs}\right\}\subset \Hom_{\C}(\C^{\alpha_i},\C^{\alpha_j})$, and $\left\{u^{ji}_{sr}\right\}\subset \Hom_{\C}(\C^{\alpha_j},\C^{\alpha_i})$,\  where $s=1,\dots,\alpha_i$ and $r=1,\dots,\alpha_j$, to be dual bases with respect to the trace pairing, that is $\tr_{\C^{\alpha_i}}(u_{rs}^{ij}u^{ji}_{lp})=\delta_{sl}\delta_{rp}$. Let us  choose the following  bases of $\Rep_{\alpha}(\overline{Q})$ compatible with the decomposition in formula (\ref{decompositionsymplectic}). 
\begin{equation}\label{goodbases}
\begin{array}{l}  
\left\{\rho_{a,rs}\right\},\mbox{\ where\ } a\in \overline{Q}, \mbox{\ and\ } \rho_{a,rs}(b)=\left\{\begin{array}{ll}\phantom{-}u_{rs}^{t(a)h(a)}&\mbox{if\ }b=a \\
\phantom{-}0&\mbox{else}\end{array}\right.\\
\left\{\rho^{a,rs}\right\},\mbox{\ where\ } a\in \overline{Q}, \mbox{\ and\ } \rho^{a,rs}(b)=\left\{\begin{array}{ll}-u_{rs}^{t(a)h(a)}&\mbox{if\ } b=a \mbox{\ and\ } a\in Q\\\phantom{-}u_{rs}^{t(a)h(a)}&\mbox{if\ } b=a \mbox{\ and\ } a\in Q^{op}\\ \phantom{-}0 &\mbox{else}\end{array}\right.
\end{array}
\end{equation}

It is easy to calculate that 
$$
\eta\left(\rho_{a,rs},\rho^{b,lp}\right)=\delta_{sl}\delta_{rp}(\delta_{(b=a^{\ast}, a\in Q)}+\delta_{(a=b^{\ast}, b\in Q)})
$$
so that the above bases, taken in the appropriate order, are dual to each other with respect to the symplectic form $\eta$.

Using these bases, we can immediately see that the operator $F$ in (\ref{operatorrelationII}), which is clearly independent on the choice of dual bases, is zero unless  $a=b^{\ast}$ and $ b\in Q$ or $b=a^{\ast}$ and $a\in Q$. 
Suppose now that $a=b^{\ast}$ and $b\in Q$. Then we have
\begin{align*}
F=\sum_q\rho_q(b^{\ast})\otimes\rho^q(b)=&\sum_{\begin{array}{c}\scriptscriptstyle{s=1,\dots,h(b)}\\\scriptscriptstyle{r=1,\dots,t(b)}\end{array}}\rho_{b^{\ast},rs}(b^{\ast})\otimes\rho^{b,sr}(b). 
\end{align*}
This map is zero on $U_i\otimes U_j$ if $i\neq t(b^{\ast})=t(a)$ and $j\neq t(b)=h(a)$. For $i=t(a)$, $j=h(a)$ we have
\begin{align*}
F=-\sum_{\begin{array}{c}\scriptscriptstyle{s=1,\dots,h(a)}\\\scriptscriptstyle{r=1,\dots,t(a)}\end{array}}u^{t(a)h(a)}_{rs}\otimes u_{sr}^{h(a)t(a)}:U_{t(a)}\otimes U_{h(a)}\longrightarrow U_{h(a)}\otimes U_{t(a)}. 
\end{align*}
By the fact that $\left\{u^{t(a)h(a)}_{rs}\right\}$ and $\left\{u_{sr}^{h(a)t(a)}\right\}$ are dual with respect to the trace pairing, the above map is $-1$ times the flip switching the two factors. Similarly, if $b=a^{\ast}$ and $a\in Q$ we get the flip, this time with positive sign.
\epf

We are now ready to prove Theorem \ref{representationtheorem}.

\textbf{Proof (of Theorem \ref{representationtheorem}).}\ 
The commutation relations between the operators $(a)_{\ell}$ and the elements of $S_n$ are easily verified. It remains to check that the commutation relations $(I)$ and $(II)$ of Definition \ref{ggsimple} hold.

For relation $(I)$ we have  
\begin{align*}
\sum_{a\in Q}[(a)_{\ell}, (a^{\ast})_{\ell}]=\sum_{p,q}\sum_{a\in Q}\rho_p|_{\M}\rho_q|_{\M}\left(\rho^p(a)\rho^q({a}^{\ast})\right)_{\ell}-\rho_q|_{\M}\rho_p|_{\M}\left(\rho^q({a}^{\ast})\rho^p(a)\right)_{\ell}.
\end{align*}
By adding and subtracting the terms 
$$
\frac{1}{2}\sum_{p,q}\sum_{a\in Q}\rho_q|_{\M}\rho_p|_{\M}\left(\rho^p(a)\rho^q(a^{\ast})\right)_{\ell}\mbox{\ and\ } \frac{1}{2}\sum_{p,q}\sum_{a\in Q}\rho_p|_{\M}\rho_q|_{\M}\left(\rho^q({a}^{\ast})\rho^p(a)\right)_{\ell}
$$
we can rewrite the above expression as
\begin{align*}
\sum_{a\in Q}[(a)_{\ell}, (a^{\ast})_{\ell}]=&\frac{1}{2}\sum_{p,q}\sum_{a\in Q}\left[\rho_p|_{\M},\rho_q|_{\M}\right]\left(\rho^p(a)\rho^q(a^{\ast})\right )_{\ell}+\left[\rho_p|_{\M},\rho_q|_{\M}\right]\left(\rho^q(a^{\ast})\rho^p(a)\right )_{\ell}\\
&\qquad\qquad\qquad+\frac{1}{2}\sum_{a\in Q}\sum_{p,q}\left(\rho_p|_{\M}\rho_q|_{\M}+\rho_q|_{\M}\rho_p|_{\M}\right)\left([\rho^p(a),\rho^q(a^{\ast})]\right)_{\ell}.
\end{align*}
Using the defining relation of $W_{\alpha}$ the first summand can be rewritten as
\begin{align*}
&\frac{1}{2}\sum_{p,q}\sum_{a\in Q}\eta(\rho_p,\rho_q)\left(\rho^p(a)\rho^q(a^{\ast})\right )_{\ell}+\eta(\rho_p,\rho_q)\left(\rho^q(a^{\ast})\rho^p(a)\right )_{\ell}\\
=&\frac{1}{2}\sum_{a\in Q}\sum_q\left(\left(\sum_p\eta(\rho_p,\rho_q)\rho^p(a)\right)\rho^q({a}^{\ast})\right)_{\ell}+\sum_p\left(\left(\sum_q\eta(\rho_p,\rho_q)\rho^q({a}^{\ast})\right)\rho^p(a)\right)_{\ell}\\
=&\frac{1}{2}\sum_{a\in Q}\sum_q\left(\rho_q(a)\rho^q({a}^{\ast})-\rho_q({a}^{\ast})\rho^q(a)\right)_{\ell}.
\end{align*}
Thus we have
\begin{align}\label{relationI}
\sum_{a\in Q}[(a)_{\ell}, (a^{\ast})_{\ell}]=&\frac{1}{2}\underbrace{\sum_{a\in Q}\sum_q\left(\rho_q(a)\rho^q({a}^{\ast})-\rho_q({a}^{\ast})\rho^q(a)\right)_{\ell}}_{A_1} \\
&\qquad\qquad\qquad+\frac{1}{2}\underbrace{\sum_{a\in Q}\sum_{p,q}\left(\rho_p|_{\M}\rho_q|_{\M}+\rho_q|_{\M}\rho_p|_{\M}\right)\left(\left[\rho^p(a),\rho^q(a^{\ast})\right]\right)_{\ell}}_{A_2}\nonumber.
\end{align}
The term $A_1$ is an operator on $U$ and it commutes with the $B$-action. This is because $A_1$  is a sum of cyclic paths, that is paths starting and ending at the same vertex, and thus commutes with each idempotent $e_i$. Moreover, $A_1$ commutes with all the elements of $\Gg(\alpha)$. To prove this  it is enough to observe that conjugating $A_1$ by an element of $g\in\Gg(\alpha)$ is the same as  changing the dual bases $\left\{\rho_q\right\}$,$\left\{\rho^q\right\}$ into their conjugate bases $\left\{g\rho_qg^{-1}\right\}$,$\left\{g\rho^qg^{-1}\right\}$. These  still form a pair of dual bases, since  the $\Gg(\alpha)$-action on $\Rep_{\alpha}(\overline{Q})$  preserves the form $\eta$. Thus this change does not affect  $A_1$, which is clearly independent on the choice of  basis. It follows that $A_1$ acts as a scalar on each irreducible component of the $B$-module $U$, and such scalar depends only on the isomorphism class of the component. Each irreducible  $B$-module is of the form $\mathcal{N}_i$ for some $i\in I$, where $\mathcal{N}_i$ is the unique irreducible of dimension vector $\epsilon_i=(0,\dots,1_i,\dots,0)$. Let us denote by $\mu_i$ the scalar by which $A_1$ acts on each irreducible component of type $\mathcal{N}_i$. To determine $\mu_i$ we can use a trace computation. Let $U_i$ be the isotypic component of $U$ of type $\mathcal{N}_i$. We know $\di U_i=\di \C^{\alpha_i}=\alpha_i$. We have
\begin{align}
\mu_i=&\frac{1}{\alpha_i}\tr|_{U_i}\left(\sum_{a\in Q}\sum_q\rho_q(a)\rho^q({a}^{\ast})-\rho_q({a}^{\ast})\rho^q(a)\right)\nonumber\\
          =&\frac{1}{\alpha_i}\sum_q\tr|_{U_i}\left(\sum_{a\in Q}\rho_q(a)\rho^q({a}^{\ast})-\rho_q({a}^{\ast})\rho^q(a)\right)\label{constanttrace}.
\end{align}
It is easy to see that we can choose dual bases $\left\{\rho_q\right\}$, $\left\{\rho^q\right\}$ compatible with the decomposition  $\Rep_{\alpha}(\overline{Q})=\bigoplus_{j,i}\Hom_{\C}\left(E_{ji},\Hom_{\C}\left(\C^{\alpha_i},\C^{\alpha_j}\right)\right)$.  Let  $D_i$ be the set of indices $q$ such that $\rho^q$ belongs to the subspace $\bigoplus_{j}\Hom_{\C}\left(E_{ji},\Hom_{\C}\left(\C^{\alpha_i},\C^{\alpha_j}\right)\right)$. It is clear that the only  terms that contribute to the  sum \eqref{constanttrace} are the ones with  $q\in D_i$ . Thus 
\begin{align*}
\mu_i=&\frac{1}{\alpha_i}\sum_{q}\tr|_{U_i}\left(\sum_{a\in Q}\rho_q(a)\rho^q({a}^{\ast})-\rho_q({a}^{\ast})\rho^q(a)\right)\\
=&\frac{1}{\alpha_i}\sum_{q\in D_i}\tr|_{U_i}\left(\sum_{a\in Q}\rho_q(a)\rho^q({a}^{\ast})-\rho_q({a}^{\ast})\rho^q(a)\right)\\
=&\frac{1}{\alpha_i}\sum_{q\in D_i}\sum_{a\in Q}\tr\left(\rho_q(a)\rho^q({a}^{\ast})-\rho^q(a)\rho_q({a}^{\ast})\right)\\
=&\frac{1}{\alpha_i}\sum_{q\in D_i}\eta(\rho_q,\rho^q)\\
=&\frac{1}{\alpha_i}\di\left(\bigoplus_{j}\Hom_{\C}\left(E_{ji},\Hom_{\C}\left(\C^{\alpha_i},\C^{\alpha_j}\right)\right)\right)\\
=&\frac{1}{\alpha_i}\sum_ja_{ji}\alpha_i\alpha_j\\
=&\sum_ja_{ji}\alpha_j.
\end{align*}
Thus for the first summand on the right hand side of \eqref{relationI} we have 
\begin{equation}\label{A1}
\frac{1}{2}A_1=\frac{1}{2}\left(\sum_{i\in I}\mu_ie_i\right)_{\ell}.
\end{equation}
Let us now look at the second summand. We have
\begin{align*}
\frac{1}{2}A_2=&\frac{1}{2}\sum_{p,q}\left(\rho_p|_{\M}\rho_q|_{\M}+\rho_q|_{\M}\rho_p|_{\M}\right)\left(\sum_{a\in Q}\left[\rho^p(a),\rho^q(a^{\ast})\right]\right)_{\ell}\\
=&\frac{1}{4}\sum_{p,q}\left(\rho_p|_{\M}\rho_q|_{\M}+\rho_q|_{\M}\rho_p|_{\M}\right)\left(\sum_{a\in Q}\left[\rho^p(a),\rho^q(a^{\ast})\right]-\left[\rho^p(a^{\ast}),\rho^q(a)\right]\right)_{\ell}
\end{align*}
The operator $\sum_{a\in Q}\left[\rho^p(a),\rho^q(a^{\ast})\right]-\left[\rho^p(a^{\ast}),\rho^q(a)\right]$ on  $U$  is a sum of cyclic paths and  it is thus in $\mathfrak{g}(\alpha)$. It follows that it can be written in the form
\begin{equation*}
\sum_{k}\tr|_U\left(\left(\sum_{a\in Q}\left[\rho^p(a),\rho^q(a^{\ast})\right]-\left[\rho^p(a^{\ast}),\rho^q(a)\right]\right)x_k\right)x^k
\end{equation*} 
where $\left\{x_k\right\}$, $\left\{x^k\right\}$ are a pair of dual bases of $\mathfrak{g}(\alpha)$ with respect to the trace form.
We can calculate
\begin{align*}
\tr|_U\!\!\left(\!\sum_{a\in Q}\left[\rho^p(a),\rho^q(a^{\ast})\right]x_k-\left[\rho^p(a^{\ast}),\rho^q(a)\right]x_k\!\!\right)\!\!=&\tr|_U\!\!\left(\!\sum_{a\in Q}\left[\rho^q(a^{\ast}),x_k\right]\rho^p(a)-\left[\rho^q(a),x_k\right]\rho^p(a^{\ast})\!\!\right)\\
=&\!-\!\sum_{a\in Q}\tr\left(\rho^p(a)(x_k\rho^q)(a^{\ast})-(x_k\rho^q)(a)\rho^p(a^{\ast})\right)\\
=&\!-\!\eta\left(\rho^p,x_k\rho^q\right).
\end{align*}
Thus we have
\begin{align}
\frac{1}{2}A_2=&-\frac{1}{4}\sum_k\sum_{p,q}\left(\rho_p|_{\M}\rho_q|_{\M}+\rho_q|_{\M}\rho_p|_{\M}\right)\left(\eta\left(\rho^p,x_k\rho^q\right)x^k\right)_{\ell}\nonumber\\
=&-\frac{1}{4}\sum_k\sum_{p,q}\left(\rho_p|_{\M}\rho_q|_{\M}\eta\left(\rho^p,x_k\rho^q\right)+\rho_q|_{\M}\rho_p|_{\M}\eta\left(\rho^p,x_k\rho^q\right)\right)(x^k)_{\ell}\nonumber\\
=&\frac{1}{4}\sum_k\sum_q\left((x_k\rho^q)\rho_q+\rho_q(x_k\rho^q)\right)|_M(x^k)_{\ell}\nonumber\\
=&\sum_k\xi(x_k)|_M(x^k)_{\ell}\nonumber\\
=&\sum_k\chi(x_k)(x^k)_{\ell}-\sum_k\sum_{m}(x_k)_m(x^k)_{\ell}\nonumber\\
=&\sum_k\chi(x_k)(x^k)_{\ell}-\sum_k(x_k)_{\ell}(x^k)_{\ell} -\sum_k\sum_{m\neq l}(x_k)_m(x^k)_{\ell}\label{A2}.
\end{align}
 We will now analyze separately the three summands in  (\ref{A2}). Observe now that, as we already mentioned, if for an element $x\in\mathfrak{g}(\alpha)$ we  write $x=\left\{x[i]\right\}_{i\in I}$ where $x[i]\in\mathfrak{gl}_{\alpha_i}$, then we have  
 \begin{equation}\label{character}
 \chi(x)=\sum_{i\in I}\chi_i\tr(x[i])
 \end{equation}
The first summand $\sum_k\chi(x_k)(x^k)_{\ell}$ commutes with all the elements of $\Gg(\alpha)$. Indeed, since $\chi$ is a character, conjugating    $\sum_k\chi(x_k)(x^k)_{\ell}$ by an element of $\Gg(\alpha)$ is the same as changing the dual bases $\left\{x_k\right\}$, $\left\{x^k\right\}$  into their conjugate bases, and this does not affect the operator. Thus, for every $i\in I$, the element $\sum_k\chi(x_k)(x^k)_{\ell}$ acts by a scalar $s_i$ on the summand $U_i$ of the $\ell$-th factor $U$ in the tensor product $U^{\otimes n}$. Using (\ref{character}) and a trace computation,  it is easy to see  that $s_i=\chi_i$. 
Similarly one can calculate that the second summand $\sum_k(x_k)_{\ell}(x^k)_{\ell}$ acts as the scalar $\alpha_i$ on $U_i$.
As for the last  summand $\sum_k(x_k)_m(x^k)_{\ell}$ it is enough to observe that the operator $\sum_kx_k\otimes x^k$ acts on  $U_i\otimes U_j$ by zero if $i\neq j$, and by the \textit{flip} switching the two factors if  $i=j$. Thus  we have that 
$\sum_k(x_k)_m(x^k)_{\ell}=P_{\ell,m}s_{\ell m}$. Adding the three summands together we get:
\begin{equation}\label{A22}
\frac{1}{2}A_2=\left(\sum_{i}(\chi_i-\alpha_i)e_i\right)_{\ell}-\sum_{m\neq \ell}P_{\ell,m}s_{\ell m}. 
\end{equation}
Plugging $(\ref{A1})$ and $(\ref{A22})$ into $(\ref{relationI})$ we get relation $(I)$.

Let us now look at relation $(II)$. We have
\begin{align}
\left[(a)_l, (b)_m\right]=&\sum_{p,q}\left[\rho_p|_{\M},\rho_q|_{\M}\right]\left(\rho^p(a)\right)_{\ell}\left(\rho^q(b)\right)_m\nonumber\\
                         =&\sum_{p,q}\eta\left(\rho_p,\rho_q\right)\left(\rho^p(a)\right)_{\ell}\left(\rho^q(b)\right)_m\nonumber\\
                         =&\sum_{q}\left(\sum_p\eta\left(\rho_p,\rho_q\right)\rho^p(a)\right)_{\ell}\left(\rho^q(b)\right)_m\nonumber\\
                         =&\sum_q\left(\rho_q(a)\right)_{\ell}\left(\rho^q(b)\right)_m\label{operatorrelationIIbis}.
\end{align}
If we regard  the operator in \eqref{operatorrelationIIbis}  as an operator on $U\otimes U$ we have $\sum_q\left(\rho_q(a)\right)_{\ell}\left(\rho^q(b)\right)_m=F$, where $F$ is as defined in formula \eqref{operatorrelationII}. It follows from Lemma \ref{operatorflip} that $(a)_{\ell}$, $(b)_{m}$ satisfy  relation $(II)$ for $\nu=-1$.
\epf
\subsection{Representations of the rational the Cherednik algebras  of type $A_{n-1}$ and $BC_n$.}$\phantom{a}$\\

\subsubsection{\bf Type $A_{n-1}$} If $Q$ is the Jordan quiver $\tilde{A}_0$ and  $\alpha=N$ we have that $\mathrm{Rep}_{N}(\tilde{A}_0)=\mathfrak{gl}_N$. In this case the conditions on the  parameters in Theorem \ref{representationtheorem} give $\chi=\lambda=n/N$ and $\nu=-1$.  We thus get a functor $F_{n,n/N}:\mathcal{D}(\mathfrak{gl}_N)-\mathrm{mod}\longrightarrow \mathcal{A}_{n,n/N,-1}(\tilde{A}_0)-\mathrm{mod}$. 

As we observed in Section \ref{sectioncherednik},  the algebra $\mathcal{A}_{n,n/N,-1}(\tilde{A}_0)$ contains a subalgebra isomorphic to the Cherednik algebra  $H_{1,N/n}(S_n,\C^{n-1})$. It follows that, by restriction,  we get a functor $F_{n,n/N}':\mathcal{D}(\mathfrak{gl}_N)-\mathrm{mod}\longrightarrow H_{1,N/n}(S_n,\C^{n-1})-\mathrm{mod}$. 

The underlying vector space of a  $H_{1,N/n}(S_n,\C^{n-1})$-module in the image of this functor  is of the form  $\left(M\otimes (\C^N)^{\otimes n}\right)^{\mathfrak{gl}_N}_{n/N}$, where $M$ is any $\mathcal{D}(\mathfrak{gl}_N)-\mathrm{mod}$. Since the center of $\mathfrak{gl}_N$ acts trivially on $M$, the $n/N$-equivariant space  $(M\otimes (\C^N)^{\otimes n})^{\mathfrak{gl}_N}_{\chi}$  coincides with the $\mathfrak{sl}_N$-invariant space $\left(M\otimes (\C^N)^{\otimes n}\right)^{\mathfrak{sl}_N}$. 
We thus have the following result.
\begin{coro}\label{cheredniktypeA} Let $n$ and $N$ be in $\N$. The construction of Theorem \ref{ggrep} yields a functor $F_{n,n/N}':\mathcal{D}(\mathfrak{gl}_N)-\mathrm{mod}\longrightarrow H_{1,N/n}(S_n,\C^{n-1})-\mathrm{mod}$ defined by the formula
$$
F_{n,n/N}'(M)=(M\otimes(\C^N)^{\otimes n})^{\mathfrak{sl}_N}.
$$
\end{coro}
\epf 

It is not hard to see that  the action of the subalgebra $H_{1,N/n}(S_n,\C^{n-1})$ of $\mathcal{A}_{n,n/N,-1}(\tilde{A}_0)$ on $F_{n,n/N}'(M)$ can be described in terms of the action of the subalgebra $\mathcal{D}(\mathfrak{sl}_N)$ of $\mathcal{D}(\mathfrak{gl}_N)$  only.   
It is now easy to check that if we regard   $M$ as a $\mathcal{D}(\mathfrak{sl}_N)$-module  the functor $F_{n,n/N}'$ agrees with the functor $F_n$ defined in  \cite{CEE}, Proposition 8.1 and Section 9.6.    
\subsubsection{\bf Type $BC_n$.} Suppose now that $Q$ is the cyclic quiver of length two. Fix a positive integer $N$ and let $(p,q)$ be a dimension vector such that $p+q=N$.  In this case we have  that $\mathfrak{gl}((p,q))=\mathfrak{gl}_p\times\mathfrak{gl}_q$ and  we can identify the algebra $\mathfrak{gl}((p,q))$ with a subalgebra of $\mathfrak{gl}_N$ in the usual way. Following the notation of \cite{EFM} Section 4.1, let us set $\mathfrak{gl}_N=\mathfrak{g}$ and let us write $\mathfrak{k}=\mathfrak{gl}((p,q))$.  We have an isomorphism 
$$
\Rep_{(p,q)}(Q)=\mathrm{Hom}(\C^p,\C^q)\oplus\mathrm{Hom}(\C^q,\C^p)\simeq \mathfrak{g}/\mathfrak{k}
$$
from which it follows that
$$
\mathcal{D}(\mathrm{Rep}_{(p,q)}(Q))\simeq\mathcal{D}(\mathfrak{g}/\mathfrak{k}).
$$
Let now $\chi=\chi_1\tr_{\mathfrak{gl}_p}+\chi_2\tr_{\mathfrak{gl}_q}$ be a character of $\mathfrak{k}$ satisfying condition \eqref{conditionchis} that is $\chi_1p+\chi_2q=n$. It is clear that such a character can always be written in the form
\begin{equation}\label{characterBCn}
\chi=\mu(q\tr_{\mathfrak{gl}_p}-p\tr_{\mathfrak{gl}_q})+\frac{n}{N}(\tr_{\mathfrak{gl}_p}+\tr_{\mathfrak{gl}_q})
\end{equation}
where $\mu\in \C$, and  $q\tr_{\mathfrak{gl}_p}-p\tr_{\mathfrak{gl}_q}$ is a basis for the subspace of  characters of $\mathfrak{k}$ that vanish on the scalars. 
Let $\mathfrak{k}_0$ be the subalgebra of traceless matrices in $\mathfrak{k}$. Clearly we have that $\mathfrak{k}=\mathfrak{k}_0\oplus \mathfrak{s}$, where $\mathfrak{s}\simeq \C$ denotes the subalgebra of scalars.  Let now $M$ be any $\mathcal{D}(\mathfrak{g}/\mathfrak{k})$-module, and let us identify the $\mathfrak{k}$-module $U=\C^p\oplus\C^q$ with $\C^N$.  Since $\mathfrak{s}$ acts trivially on $M$ we have that
$$
(M\otimes(\C^N)^{\otimes n})^{\mathfrak{k}}_{\chi}=(M\otimes(\C^N)^{\otimes n})^{\mathfrak{k}_0}_{\mu},
$$ 
where  $(M\otimes(\C^N)^{\otimes n})^{\mathfrak{k}_0}_{\mu}$  denotes the equivariant subspace under the action of $\mathfrak{k}_0$ with respect to the character $\mu(q\tr_{\mathfrak{gl}_p}-p\tr_{\mathfrak{gl}_q})$. 
Observe now that with the notation of formulas \eqref{characterBCn} we must have
$$
\chi_1=\mu q+\frac{n}{N},\qquad\qquad\chi_2=-\mu p+\frac{n}{N}.
$$
Moreover we have that for this quiver the matrix of the symmetrized Ringel form is
$$
C=\left(\begin{array}{cc} 2 &-2\\
                   -2 & 2\end{array}\right).
$$ 
It follows that in this case the values of the parameters $\lambda$ and $\nu$ for which the representation of Theorem \ref{representationtheorem} exists are given by the following formulas
\begin{equation}
\lambda_1=\mu q + \frac{n}{N} -p+q,\qquad\qquad
\lambda_2=-\mu p+\frac{n}{N}+p-q,\qquad\qquad
\nu=-1.
\end{equation}
Recall now (Section \ref{sectioncherednik}) that we have an isomorphism $\mathcal{A}_{n,\lambda,\nu}(Q)\simeq H_{t,k,c}((\Z/2\Z)^n\rtimes S_n, \C^n)$ where the parameters $(\lambda,\nu)$ and $(t,k,c)$ are related by the formulas in \eqref{parametersBCn} so that
\begin{equation}
t=\frac{n}{N}+\frac{\mu(q-p)}{2},\qquad\qquad k=1,\qquad\qquad c=\frac{\mu N}{2}+(q-p).\nonumber
\end{equation}
Rescaling the parameters by a factor $2$ we have the following result.
\begin{coro} Let $p$ and $q$ be in $\N$ with $p+q=N$ and let $\mu$ be in $\C$. The construction of Theorem \ref{ggrep} yields  a functor $F_{n,p,\mu}: \mathcal{D}(\mathfrak{g}/\mathfrak{k})-\mathrm{mod}\rightarrow H_{t,k,c}((\Z/2\Z)^n\rtimes S_n, \C^n)-\mathrm{mod}$ given by the formula
$$
F_{n,p,\mu}(M)=(M\otimes(\C^N)^{\otimes n})^{\mathfrak{k}_0}_{\mu}
$$
where the values of the parameters  are as follows
\begin{equation}\label{parametersBCncomparison}
t=\frac{2n}{N}+\mu(q-p),\qquad\qquad k=2,\qquad\qquad c=\mu N+2(q-p).
\end{equation}
\end{coro}
\epf

Let us now set $G=\mathrm{GL}_N$ and $K=\mathrm{GL}_p\times \mathrm{GL}_q$ so that $\mathfrak{g}=\mathrm{Lie}(G)$ and $\mathfrak{k}=\mathrm{Lie}(K)$. Let $\lambda$ be a complex number and let $\lambda(q\tr_{\mathfrak{gl}_p}-p\tr_{\mathfrak{gl}_q})$ be a character of $\mathfrak{k}$. Recall that in Section 5.1 of  \cite{EFM}  the authors define a functor $F_{n,p,\mu}^{\lambda}$ from the category of modules over the algebra of twisted differential operators $\mathcal{D}^{\lambda}(G/K)$   to the category of modules over the dDAHA $\HH(k_1,k_2,k_3)$ with parameters
$$
t=\frac{2n}{N}+(\lambda+\mu)(q-p),\qquad k_1=1,\qquad k_2=p-q-\lambda N,\qquad k_3=(\lambda-\mu)N.
$$
Taking into account  the relations between the parameters $(k,c)$ and $(k_1,k_2,k_3)$ in formulas \eqref{parametersdegeneration}, it is not hard to see that the functor $F_{n,p,\mu}$ can be obtained from the functor $F^{\lambda}_{n,p,\mu}$ by rational degeneration when $\lambda=0$.

\section{Rational generalized double affine Hecke algebras of higher rank}\label{gdahadefisec}

\subsection{Definition of the rational GDAHA of higher rank.} Let us recall the definition of the rational generalized double affine Hecke algebra (rational GDAHA) of higher rank introduced in \cite{EGO}. Let $m$ be a positive integer. Let $D$ be an $m$-branched star-shaped graph, that is a tree with one $m$-valent vertex, called the node or branching vertex and labeled by $0$ , and the rest of the vertices $2$- and $1$-valent.  We will assume that $D$ is \emph{not} a finite Dynkin diagram. The graph $D$ has $m$ branches (legs); denote their lengths plus one by $\ell_1,\dots,\ell_m$. For   $k=1,\dots,m$ and $j=1,\dots,\ell_k$, let $\g_{kj}$ be complex parameters, and let $\g$ denote the collection of all parameters $\g_{kj}$. Let $\nu$ be an additional complex parameter. 
\begin{defi}\label{gdahadefi}
The rational (or degenerate) generalized DAHA of rank $n$ attached to $D$ and to the parameters $\g, \nu$ is the algebra $B_n(\gamma,\nu)$ generated over $\C$ by elements $Y_{i,k}$ (where $i=1,\dots,n$ and $k=1,\dots,m$) and the symmetric group $S_n$, with  defining relations (r1)-(r6) as follows. For any $i,j,h\in[1,n]$ with $i\neq j$, and any $k,l\in [1,m]$:
\begin{align*}
(r1)\qquad & s_{ij}Y_{i,k}=Y_{j,k}s_{ij},\\
(r2)\qquad & s_{ij}Y_{h,k}=Y_{h,k}s_{ij},\;if \; h\neq i,j,\\
(r3)\qquad &\prod_{j=1}^{\ell_k}(Y_{i,k}-\g_{kj})=0,\\ 
(r4)\qquad & Y_{i,1}+Y_{i,2}+\dots+Y_{i,m}=\nu\sum_{j\neq i}s_{ij},\\
(r5)\qquad & [Y_{i,k},Y_{j,k}]=\nu(Y_{i,k}-Y_{j,k})s_{ij},\\
(r6)\qquad & [Y_{i,k},Y_{j,l}]=0,\; if \;k\neq l 
\end{align*}
\end{defi}

Let us  label the vertices of $D$, excluding the node, by pairs $(k,s)$ where $ k=1,\dots, m$ is the leg-number, and $s=1,\dots, \ell_k-1$ denotes the position of the vertex on the leg  starting from the inside  (so that $(k, 1)$ is connected to the node).  As explained in \cite[$\S\; 2.1$]{EGO}, there exist unique complex parameters $\lambda_{0}=\lambda_{0}(\g)$, $\lambda_{(k,s)}=\lambda_{(k,s)}(\g)$,  and  $\xi_{k}=\xi_k(\g)$  such that $\xi_1+\dots+\xi_m=0$ and the equation
\begin{equation}\label{parametersgdaha1}
\g_{kj}=\sum_{s=1}^{j-1}\lambda_{(k,s)}+\frac{\lambda_{0}}{m}+\xi_k.
\end{equation}
is satisfied for all $j,k$.
We  denote by $\lambda$ the collection of all parameters $\lambda_{0},\lambda_{(k,s)}$.  As explained in \cite{EGO},  for any choice of complex parameters $\sigma_1,\dots,\sigma_m$ such that $\sigma_1+\dots +\sigma_m=0$ the  transformation $\g_{kj}\rightarrow \g_{kj}+\sigma_k$ leaves the algebra  $B_n(\g,\nu)$ unchanged up to isomorphism  (the required isomorphism is given by the assignment $Y_{i,k}\rightarrow Y_{i,k}+\sigma_k$). Thus the parameters  $\xi_1,\dots,\xi_m$ are ``redundant'';  their use though allows to simplify the presentation.
\subsection{Relation to the Gan-Ginzburg algebra.}
Let now $Q$ be the quiver obtained from $D$ by giving to all edges the orientation pointing toward the node. For each $ k=1,\dots,  m$ and $ s=1,\dots, \ell_k-1$ denote by $h_{k, s}$ the arrow in the $k$-th leg with tail at the vertex $(k, s)$ and head at the vertex $(k, s-1)$, where we set $(k, 0):=0$ for all $k$. Thus in particular the path $h_{k,1}h_{k,1}^{\ast}$ is the  unique path of length two along the $k$-th leg that starts and ends at the node. Let $e_{0}$ be the idempotent in $\C\overline Q$ corresponding to the node, and let the parameters $\gamma$ and $\lambda$ be related by (\ref{parametersgdaha1}).   We have the following result.
\begin{prop}[\cite{EGO} Proposition 2.3.4]\label{isomorphism}  There exists a natural algebra isomorphism \\$\varphi: B_n(\gamma,\nu)\simeq ~ e_{0}^{\otimes n}\mathcal{A}_{n, \lambda, -\nu}(Q)e_{0}^{\otimes n}$ defined by the formulas
$$
\varphi(Y_{i,k})=e_{0}^{\otimes i-1}\otimes\left(-h_{k1}h_{k1}^{\ast}+\left(\xi_k+\frac{\lambda_{0}}{m}\right)e_{0}\right)\otimes e_{0}^{\otimes n-i}.
$$
\end{prop}

\epf
\begin{rem}\label{remarksign}
In  Proposition 2.3.4 of \cite{EGO} there is a minor sign mistake, and the above isomorphism figures as $B_n(\gamma,\nu)\simeq e_{0}^{\otimes n}\mathcal{A}_{n, \lambda, \nu}(Q)e_{0}^{\otimes n}$. The proof of the proposition stays the same once the sign change in the parameter $\nu$ is considered.  
\end{rem}
\begin{rem}\label{remarkorientation}
The presentation of the algebra $B_{n}(\nu,\g)$ depends only on the graph $D$ and not on the orientation of its edges.  As we observed before,  the algebra $A_{n, \lambda,\nu}(Q)$  does not depend up to isomorphism on the orientation of the edges of $Q$, but only on the underlying graph: if the orientation of an arrow $a\in Q$ is changed then the required isomorphism is given by  $(a)_l\rightarrow(a^{\ast})_l$, $(a^{\ast})_{\ell}\rightarrow -(a)_{\ell}$ for $\ell\in[1,n]$. On the contrary,  the presentation of $A_{n, \lambda, \nu}(Q)$ by generators and relations  \emph{does} depend on the orientation. That is why, in order to give explicit formulas for the isomorphism $\varphi$,  we assigned an orientation to the graph $D$. Nevertheless one can say that the isomorphism $\varphi$  is ``natural'':  if the orientation of the  edge connecting the $k$-th leg to the node  is changed, then the sign with which the corresponding length-two loop at the node occurs in the defining relations changes, and $\varphi$ is modified accordingly by substituting $-h_{k,1}h_{k,1}^{\ast}$ with $h_{k,1}^{\ast}h_{k,1}$ in the formula of Proposition \ref{isomorphism}.   \end{rem}
\section{ Representations of  rational GDAHA}\label{gdahasrep}
In this section we will give an explicit Lie-theoretic construction of some finite dimensional representations for rational GDAHA of higher rank. 

For a positive integer $N$, let $\left\{\e_i\right\}_{i=1,\dots,N}$ be the standard basis of the weight lattice $\Z^N$ of  $\gln$. Then the $i$-th fundamental weight for $\gln$ is $\omega_i=\e_1+\dots+\e_i$. We  say that a finite dimensional irreducible representation of $\mathfrak{gl}_N$ is $\ell$-stepped if its highest weight $\mu$ is of the form $\mu=a_1\omega_{p_1}+\dots+a_{\ell}\omega_{p_{\ell}}$, where $p_{\ell}=N$. In other words the highest weight contains at most $\ell$ distinct entries.  Recall that if $V$ is an $\ell$-stepped $\gln$-module, and $\C^N$ is the vector  representation of $\gln$, then $V\otimes \C^N$ is a direct sum of at most $\ell$ distinct  irreducible representations.

Let  $V_1,\dots, V_m$ be irreducible finite dimensional representations  of $\gln$ such that $V_k$ is $\ell_k$~-~stepped but not $(\ell_k-1)$~-~stepped, and let $\C^N$ be the vector representation of this Lie algebra. Let $\chi'$ be a character of $\gln$. We define the $\chi'$-equivariant subspace of the tensor product  $V_1\otimes\cdots\otimes V_m\otimes (\C^N)^{\otimes n}$ as 
$$
E_{n,\chi'}(V_1,\dots, V_m):=\left\{v\in V_1\otimes\cdots\otimes V_m\otimes (\C^N)^{\otimes n}  | xv=\chi'(x)v,\,\;\forall\, x\in \gln  \right\}.
$$ 

Let $\Omega=\sum \U_{rs}\otimes \U_{sr}$  be the Casimir tensor of $\gln$, where $\U_{rs}$ denotes the corresponding matrix unit in $\mathrm{Mat}_{N\times N}(\C)$. Let $\Omega_{ij}$ denote the tensor $\Omega$ acting in the $i$-th and $j$-th factors of the above tensor product. Observe that if $i$ and $j$ are both greater than $m$, then $\Omega_{ij}$ is just the transposition operator. Let $c$ be a complex number such that $\chi'=c\tr_{\gln}$. Let $\nu$ be as in Section \ref{gdahadefisec}. The following theorem holds.
\begin{thm}\label{repreBn}
For an appropriate value of the parameter $\gamma$ the formulas 
\begin{equation}
s_{ij}=\Omega_{m+i,m+j} 
\end{equation}
\begin{equation}\label{repreBnformula}
Y_{i,k}=-\nu\Omega_{k,i+m}-\nu\frac{(N-c)}{m}.
\end{equation}
define a representation of the algebra $B_n(\nu,\g)$ on the space $E_{n,\chi'}(V_1,\dots,V_m)$.
\end{thm}
\pf\  Relations (r1), (r2) and (r6) are obviously satisfied for all choices of $\g$.  Relation (r3) follows from the fact that the representation $V_k$ is $\ell_k$-stepped. In this case, in fact, the module $V_k\otimes \C^N$ has  at most $\ell_k$ irreducible components, and hence the Casimir tensor $\Omega$ has at most $\ell_k$ distinct eigenvalues (and is semisimple) on this tensor product (as it acts by a scalar on each summand).  Relation (r4) follows from the fact that we work in the $\chi'$-equivariant subspace, and the fact that the standard quadratic Casimir of $\gln$ has eigenvalue $N$ on the vector representation.  Indeed if for $i=1,\dots,n$ we denote  by $C|_{i+m}$ the quadratic Casimir for $\gln$ acting in the $i+m$-th factor of the tensor product we have that  $C|_{i+m}=N$. Thus we have
\begin{align*}
Y_{i,1}+\dots+Y_{i,m}|_{{\scriptscriptstyle E_{n,\chi'}(V_1,\dots,V_m)}}=&-\nu\sum_{k= 1}^m\Omega_{k,i+m}-\nu (N-c)\\
=&-\nu\Omega_{1,i+m}-\nu\sum_{k=2}^m\Omega_{k,i+m}-\nu (N-c)\\
=&\nu\sum_{j\neq i}\Omega_{j+m,i+m}+\nu C|_{i+m}+\nu\sum_{k=2}^m\Omega_{k,i+m}-\nu c+\\
&\qquad\qquad\qquad\qquad\qquad\qquad\qquad\qquad-\nu\sum_{k= 2}^m\Omega_{k,i+m}-\nu N+\nu c\\
=&\nu\sum_{j\neq i}\Omega_{j+m,i+m}+\nu N-\nu N\\
=&\nu\sum_{j\neq i}\Omega_{j+m,i+m}\\
=&\nu\sum_{j\neq i}s_{ij}|_{{\scriptscriptstyle E_{n,\chi'}(V_1,\dots,V_m)}}.
\end{align*} 
Relation (r5) can be rewritten as 
\begin{equation}\label{r5bis}
[Y_{i,k}, Y_{j,k}-\nu s_{ij}]=0.
\end{equation}
This follows from the fact that, using (r1) we can write
\begin{align*}
[Y_{i,k}, Y_{j,k}]-\nu(Y_{i,k}-Y_{j,k})s_{ij}&=[Y_{i,k}, Y_{j,k}]-\nu Y_{i,k}s_{ij}+\nu s_{ij} Y_{i,k}\\
&=[Y_{i,k}, Y_{j,k}-\nu s_{ij}].
\end{align*}
Since we have that
\begin{align*}
[Y_{i,k}, Y_{j,k}-\nu s_{ij}]_{{\scriptscriptstyle E_{n,\chi'}(V_1,\dots,V_m)}}&=\nu[\Omega_{k,i+m}, \Omega_{k,j+m}+\Omega_{i+m,j+m}]
\end{align*}
relation \eqref{r5bis} follows from the invariance property of the Casimir tensor.
\epf

Now we would like to calculate the values of the parameters $\g_{kj}$ in  Theorem \ref{repreBn} or, in other words,  the eigenvalues of $\Omega$ on the irreducible components of $V_k\otimes \C^N$ for each $k$. 

Let $V$ be an $\ell$-stepped (but not $(\ell-1)$-stepped) representation of $\gln$ of highest weight
$$
\mu:=b_1\omega_{p_1}+\dots+b_{\ell}\omega_{p_{\ell}}
$$
where $p_{\ell}=N$.  Then the tensor product $V\otimes \C^N$ has a decomposition
$$
V\otimes \C^N=Y_1+\dots +Y_{\ell}
$$
where $Y_j$ is the irreducible representation of highest weight 
$$
\mu+\e_{p_{j-1}+1},
$$
with $p_0=0$. Since $\Omega|_{V\otimes\C^N}=\frac{1}{2}(C|_{V\otimes\C^N}-C|_{V}-C|_{\C^N})$, where $C$ denotes the quadratic Casimir for $\gln$, it follows that the eigenvalue of $\Omega$ on $Y_j$ is
$$
\Omega|_{Y_j}=\frac{1}{2}((\mu+\e_{p_{j-1}+1},\mu+\e_{p_{j-1}+1}+2\rho)-(\mu,\mu+2\rho)-N)=b_j+\dots+b_{\ell}-p_{j-1},
$$
where as usual $\rho$ is half the  sum of the positive roots.

Thus, if the highest weight of $V_k$ is
$$
\mu_k:=b_{k1}\omega_{p_{k1}}+\dots+b_{k\ell_k}\omega_{p_{k\ell_k}},
$$
the above computation proves the following result.
\begin{lem}\label{Bneigenvalues} The value of the parameter $\g$ for which the representation in Theorem \ref{repreBn} exists is given by the formulas  
$$
\g_{kj}=-\nu(b_{kj}+\dots+b_{k\ell_k}-p_{k\,j-1}+(N-c)/m),
$$
where $k=1,\dots,m$ and $j=1,\dots,\ell_k$.
\end{lem}
\epf

Lemma \ref{Bneigenvalues} implies that the parameters $\lambda_{0}$, $\lambda_{(k,s)}$, and $\xi_k$ defined by the relations in (\ref{parametersgdaha1}) are given by the formulas
\begin{align}\label{lambda0}
\lambda_{0}& =-\nu\sum_{k=1}^m\sum_{j=1}^{\ell_k}b_{kj}-\nu (N-c)\\
\label{lambdas}
\lambda_{(k,s)}&=\nu \left(b_{ks}+p_{ks}-p_{k\,s-1}\right) \\
\label{xik}
\xi_k &=\nu\left(\frac{1}{m}\sum_{h\neq k}\sum_{j=1}^{\ell_h} b_{hj}+\frac{1-m}{m}\sum_{j=1}^{\ell_k}b_{kj}\right)
\end{align}
where $k=1,\dots, m$ and $s=1,\dots,\ell_k-1$.

We conclude  by observing that in order for the space $E_{n,\chi'}(V_1,\dots,V_m)$ to be nonzero, we need the center of $\gln$ to act  on the tensor product by the character $\chi'$. This gives the condition
\begin{equation}\label{traceconditionBn}
\sum_{k=1}^m\sum_{s=1}^{\ell_k}b_{ks}p_{ks}+n-cN=0.
\end{equation}

\section{Quantum Hamiltonian reduction and  differential operators for quivers of type $A$}\label{typeAqhr}

In the last section of this paper  we will show how to recover some of the representations of rational GDAHA described in the previous section  as a special case of the representations of Gan-Ginzburg algebras constructed in Section \ref{ggrep}. In order to to this we will have to investigate  the relationship between $W_{\alpha}$-modules for a star-shaped quiver  and $\gln$-modules. This will require a preliminary  discussion of the connection between partial flag varieties and spaces of quiver representations. In this analysis, our main technical tool will be a construction called \emph{quantum Hamiltonian Reduction}, which is the   quantum analog of the Hamiltonian reduction procedure. 

\subsection{ Quantum Hamiltonian reduction and twisted differential operators.}\label{qhr}
 Let us recall the definition of the quantum Hamiltonian reduction. For more details about this construction and its properties  we refer to \cite{ELOR}, Section 1.1, and \cite{E}, Chapter 4.

Let $A$ be an associative algebra, $\mathfrak{g}$ be a Lie algebra and $\mu:\mathfrak{g}\rightarrow A$ a homomorphism of Lie algebras.
\begin{defi}\label{qhrdefi} The {\bf quantum Hamiltonian reduction} of $A$ with respect to $\mathfrak{g}$ with quantum moment map $\mu$ is   the associative  algebra 
$$
\mathcal{U}(A,\mathfrak{g},\mu)=(A/A\mu(\mathfrak{g}))^{\mathfrak{g}},
$$
where the invariants are taken with respect to the adjoint action of $\mathfrak{g}$ on $A$.
\end{defi}
An important example of quantum Hamiltonian reduction will be for us the algebra of twisted differential operators on a partial flag variety. Let us recall this construction. Let $\alpha_0,\dots,\alpha_{\ell-1}$ be a collection of positive integers satisfying the condition
$$
N=\alpha_0>\alpha_1>\dots>\alpha_{\ell-1}>0
$$
Let $\mathcal{F}\ell(\alpha_1,\dots,\alpha_{\ell-1};N)$ be the space of partial flags of type $(\alpha_1,\dots,\alpha_{\ell-1};N)$, that is   the configuration space of $\ell-1$ subspaces in $\C^N$
$$
V_{\ell-1}\subset V_{\ell-2}\subset\cdots\subset V_1 \subset V_0=\C^N
$$ 
such that $\mathrm{dim}(V_i)=\alpha_i$. Let $L_i$ be the line bundle with fiber $\wedge^{\alpha_i}(V_i)$ on $\mathcal{F}\ell(\alpha_1,\dots,\alpha_{\ell-1};N)$. Let $E$ be the total space of the principal $(\C^{\times})^{\ell-1}$-bundle corresponding to the direct sum of line bundles $L_1\oplus\dots\oplus L_{\ell-1}$. Denote  by $\mathcal{E}_i$ the Euler vector field on $E$ along the $i$-th factor of the fiber. Let $\mathrm{Vect}(E)$ denote the algebra of vector fields on $E$,  and denote by $\mathcal{E}:\C^{\ell-1}\rightarrow \mathrm{Vect}(E)$ the map that sends the standard basis vectors to the $\mathcal{E}_i$s. Let   $\tilde{\chi}=(\tilde{\chi}_1,\dots,\tilde{\chi}_{\ell-1})$ be a  vector in $\C^{\ell-1}$. 
\begin{defi}\label{twistdiffopflag} The algebra of twisted differential operators on  $\mathcal{F}\ell(\alpha_1,\dots,\alpha_{\ell-1};N)$ is the algebra
$$
\mathcal{D}_{\tilde{\chi}}(\mathcal{F}\ell(\alpha_1,\dots,\alpha_{\ell-1};N))=\mathcal{U}(\mathcal{D}(E),\C^{\ell-1}, \mathcal{E}-\tilde{\chi}).
$$
\end{defi}
When $\tilde{\chi}=0$ this construction  gives the algebra   $\mathcal{D}(\mathcal{F}\ell(\alpha_1,\dots,\alpha_{\ell-1};N))$  of differential operators on the partial flag variety. 
\subsection{ Quiver-related partial flag varieties.}\label{flagsection}
Let  $Q$ be a quiver of type $A$  with vertices labeled by the integers $0,1, \dots, {\ell-1}$, and such that all edges are oriented toward the vertex $0$. Let $N$ be a positive integer, and let  $\alpha$ be a dimension vector for $Q$ such that $\alpha_{0}=N$, and the dimensions at the remaining vertices are all non-zero and  strictly decreasing  as one moves away from the vertex $0$. In other words  
\begin{equation}\label{flagsequence}
\alpha_0=N>\alpha_1>\cdots>\alpha_{\ell-1}>0 .
\end{equation}
Observe now that to any dimension vector as above we can attach the corresponding variety $\mathcal{F}\ell( \alpha_1,\dots, \alpha_{\ell-1}; N)$ of partial flags in $\C^N$. 

As in Section \ref{ggaldefi}, let us consider the  Lie group $\Gg(\alpha)$  and its Lie algebra $\mathfrak{g}(\alpha)$ defined as follows
$$
\Gg(\alpha)=\prod_{i=0}^{\ell-1} \GL_{\alpha_i}, \qquad \mathfrak{g}(\alpha)=\prod_{i=0}^{\ell-1}\mathfrak{gl}_{\alpha_i}.
$$
Consider now the  Lie group $\Gg(\alpha)'$ and its Lie algebra $\mathfrak{g}(\alpha)'$ defined as follows
$$
\Gg(\alpha)'=\prod_{i=1}^{\ell-1}\GL_{\alpha_i}, \qquad \mathfrak{g}(\alpha)'=\prod_{i=1}^{\ell-1} \mathfrak{gl}_{\alpha_i}.
$$
We have
\begin{equation}\label{flagroupgeneral}
\Gg(\alpha)=\GL_N\times \Gg(\alpha)', \qquad \qquad \mathfrak{g}(\alpha)=\gln \times\mathfrak{g}(\alpha)'.
\end{equation}

Denote by $\mathrm{Incl}(\C^{\alpha_i}, \C^{\alpha_{i-1}})\subset\Hom(\C^{\alpha_i}, \C^{\alpha_{i-1}})$ the open subset  consisting of inclusions, and by  $X=\oplus_{i=1}^{\ell-1}\mathrm{Incl}(\C^{\alpha_i}, \C^{\alpha_{i-1}})$ the open subset of $\Rep_{\alpha}(Q)$ consisting of the representations of $Q$ for which all arrows are inclusions. 
As observed in \cite{ELOR}, the group $\Gg(\alpha)'$ acts freely on $X$, and the quotient is the flag variety $\mathcal{F}\ell(\alpha_1,\dots,\alpha_{\ell-1}; N)$. Moreover,  the closed set of non-inclusions $\Rep_{\alpha}(Q)\smallsetminus X$ has codimension greater or equal to two in $\Rep_{\alpha}(Q)$ (see \cite{ELOR}, Section 2.4).   It follows  that the algebras of differential operators over $\Rep_{\alpha}(Q)$ and $X$ are isomorphic, that is we have
$$
W_{\alpha}=\mathcal{D}(\Rep_{\alpha}(Q))\cong \mathcal{D}(X).
$$ 
Note  that the action of the remaining $\GL_N$  on $\Rep_{\alpha}(Q)$  commutes with the action of $\Gg(\alpha)'$.  

Let $\mathfrak{J}$ be the natural representation  of $G(\alpha)$ on $\Rep_{\alpha}(Q)$ by conjugation. Differentiating $\mathfrak{J}$ gives a homomorphism $d\mathfrak{J}:\mathfrak{g}(\alpha)\rightarrow W_{\alpha}$. Let $\tilde{\chi}=\sum_{i=1}^{\ell}\chi_i\tr_{\mathfrak{gl}_{\alpha_i}}$ be a character of $\mathfrak{g}(\alpha)'$. The following theorem was proved in \cite{ELOR}.
\begin{thm}\label{elortheorem}\cite[Theorem 2.4.1]{ELOR}
There is a $\GL_N$-equivariant isomorphism of algebras
$$
\mathcal{U}\left(W_{\alpha},\mathfrak{g}(\alpha)' , d\mathfrak{J}-\tilde{\chi}\right)\simeq \mathcal{D}_{\tilde{\chi}}(\mathcal{F}\ell(\alpha_1,\dots,\alpha_{\ell-1};N)).
$$
\end{thm} 

\epf

\subsection{ A quiver version of the Borel-Weil construction.}\label{borelweil}
Let $Q$ and $\alpha$   be as in the previous section. In this section we will use  Theorem \ref{elortheorem} to establish a relationship between  $W_{\alpha}$-modules and finite dimensional $\gln$-modules. Our goal is to compare the representations of Section \ref{ggrep} and Section \ref{gdahasrep}. Note that in in Section $\ref{ggrep}$ we used the homomorphism $\xi$ of the Weil representation  to embed $\mathfrak{g}(\alpha)$ into $W_{\alpha}$, while in the previous section we used   the differential $d\mathfrak{J}$ of the natural action of $\mathfrak{g}(\alpha)$ on $\Rep_{\alpha}(Q)$. Because of this, we will need to modify   Theorem \ref{elortheorem} to adapt it to our situation. 

First of all, we want to  compare  the homomorphisms $\xi$ and $d\mathfrak{J}$ by finding explicit formulas for them. For this we need a convenient description of $W_{\alpha}$. Observe that, for $Q$ and $\alpha$ as above, we have  
$$
\Rep_{\alpha}(Q)=\bigoplus_{p=1}^{\ell-1}\Hom(\C^{\alpha_p}, \C^{\alpha_{p-1}}).
$$
For $i=1,\dots,\ell-1$ denote by $h_i\in Q$ the arrow starting at the vertex $i$ and ending at the vertex ${i-1}$, and by $h_i^{\ast}\in Q^{op}$ the opposite arrow. Then the algebra of polynomial functions on $\Rep_{\alpha}(Q)$ is
$$
\C[\Rep_{\alpha}(Q)]=\C[t^{h_i}_{rs}: i=1,\dots,\ell-1,\, r=1,\dots,\alpha_{i-1},\,  s=1,\dots,\alpha_i],
$$
and the algebra of differential operators with constant coefficients on $\Rep_{\alpha}(Q)$ is
$$
\mathcal{D}_c(\Rep_{\alpha}(Q))=\C[\partial/\partial t_{rs}^{h_i}: i=1,\dots,\ell-1,\, r=1,\dots,\alpha_{i-1},\, s=1,\dots,\alpha_i].
$$
With the same notation as in (\ref{goodbases}), consider the basis for the vector space $\Rep_{\alpha}(Q)$ consisting of the vectors 
\begin{equation}\label{goodbasis1}
\rho_{h_i, rs}\qquad i=1,\dots,\ell-1,\quad r=1,\dots,\alpha_{i-1},\quad  s=1,\dots,\alpha_i
\end{equation}and the basis for the vector space $\Rep_{\alpha}(Q^{op})$ consisting  of the vectors 
\begin{equation}\label{goodbasis2}
\rho^{h^{\ast}_{i}, sr} \qquad i=1,\dots,\ell-1,\quad r=1,\dots,\alpha_{i-1},\quad s=1,\dots,\alpha_i.
\end{equation} 
Let $S(\Rep_{\alpha}(Q))$ and $S(\Rep_{\alpha}(Q^{op}))$ denote the symmetric algebras of the respective vector spaces. Then the assignments
\begin{equation}\label{isooperator}
\rho_{h_i, rs}\rightarrow \partial/\partial t^{h_i}_{rs}
\end{equation}
\begin{equation}\label{isofunction}
\rho^{h_i^{\ast}, sr}\rightarrow t^{h_i}_{rs}
\end{equation}
define  algebra isomorphisms $S(\Rep_{\alpha}(Q))\simeq \mathcal{D}_c(\Rep_{\alpha}(Q))$ and   $S(\Rep_{\alpha}(Q^{op}))\simeq \C[\Rep_{\alpha}(Q)]$ respectively. It is clear that $W_{\alpha}=S(\Rep_{\alpha}(Q))\otimes S(\Rep_{\alpha}(Q^{op}))$ as a vector space, and that the above assignments give an explicit isomorphism  $W_{\alpha}\simeq \mathcal{D}(\Rep_{\alpha}(Q))$, where $\mathcal{D}(\Rep_{\alpha}(Q))$ is the algebra of differential operators with polynomial coefficients on $\Rep_{\alpha}(Q)$.
\begin{rem}\label{Waction}
Observe that with the above identifications the natural action of $D(\Rep_{\alpha}(Q))$ on $\C[\Rep_{\alpha}(Q)]$ by differential operators correspond to the action of $W_{\alpha}$ on $S(\Rep_{\alpha}(Q^{op}))$ in which $\rho^{h_i^{\ast}, sr}$ acts by multiplication and $\rho_{h_i, rs}$ by the contraction map $\eta(\rho_{h_i, rs}, -)$. 
\end{rem}
The following lemma gives  explicit formulas for the homomorphisms $d\mathfrak{J}$ and  $\xi$.
\begin{lem}\label{xidJ}
Let $\U^{i}_{rs}$ be the $rs$-th matrix unit in $\mathrm{Mat}_{\alpha_i\times \alpha_i}(\C)=\mathfrak{gl}_{\alpha_p}\subset \mathfrak{g}(\alpha)$. Then  
\begin{align}\label{xiformula}
\xi(\U^{i}_{rs})&=
\sum_{l=1,\dots, \alpha_{i-1}}\rho^{h_i^{\ast}, rl}\rho_{h_i, ls}-
\sum_{l=1,\dots, \alpha_{i+1}}\rho^{h_{i+1}^{\ast}, ls}\rho_{h_{i+1}, rl}+\delta_{rs}\frac{1}{2}(\alpha_{i-1}-\alpha_{i+1})
\end{align}
\begin{align}\label{dJformula}
d\mathfrak{J}(E^{i}_{rs})&=\sum_{l=1,\dots, \alpha_{i-1}}\rho^{h_i^{\ast}, rl}\rho_{h_i, ls}-
\sum_{l=1,\dots, \alpha_{i+1}}\rho^{h_{i+1}^{\ast}, ls}\rho_{h_{{i+1}}, rl}
\end{align}
\end{lem}
\pf\ The proof is by direct computation. For formula (\ref{xiformula}) it is enough to apply the definition of $\xi$ given in formula (\ref{defixi}). For formula (\ref{dJformula}) one can easily differentiate the action of $G(\alpha)$ on $S(\Rep_{\alpha}(Q^{op}))$ (see for example  
\cite{H}, Lemma 3.1).
\epf

Let now $\chi=\sum_{i=1}^{\ell-1}\chi_i\tr_{\mathfrak{gl}_{\alpha_i}}$ be a character of $\mathfrak{g}(\alpha)'$, and let  $\tilde{\chi}$ be the character  given by the formula
\begin{equation}\label{chitilde2}
\tilde{\chi}:=\chi-\frac{1}{2}\sum_{i=1}^{\ell-1}(\alpha_{i-1}-\alpha_{i+1})\tr_{\mathfrak{gl}_{\alpha_i}}
\end{equation}

We have the following  modified version of  Theorem 2.4.1. in \cite{ELOR}.
\begin{thm}\label{hamiltonianred}
There is a $\GL_N$-equivariant isomorphism of algebras
$$
\mathcal{U}\left(W_{\alpha},\mathfrak{g}(\alpha)' , \xi-\chi\right)\simeq \mathcal{D}_{\tilde{\chi}}(\mathcal{F}\ell(\alpha_1,\dots,\alpha_{\ell-1};N)).
$$
\end{thm}
\pf\ By Theorem 2.4.1 in \cite{ELOR} there is a $\GL_N$-equivariant algebra isomorphism
$$
\mathcal{U}\left(W_{\alpha}, \mathfrak{g}(\alpha)' ,d\mathfrak{J}-\tilde{\chi}\right)\simeq \mathcal{D}_{\tilde{\chi}}(\mathcal{F}\ell(\alpha_1,\dots,\alpha_{\ell-1};N)).
$$
By Lemma \ref{xidJ}  we have that $\xi-\chi=d\mathfrak{J}-\tilde{\chi}$, and the theorem follows. 
\epf

 Let now $M$ be a $W_{\alpha}$-module. Then $M$ can be given a structure of $\mathfrak{g}(\alpha)$-module in two ways, via the map $\xi$ and $d\mathfrak{J}$ respectively. The above lemma  allows us to compare these two $\mathfrak{g}(\alpha)$-module structures. We  have the following proposition. 
\begin{prop}\label{mapsdifference}
For every  $W_{\alpha}$-module $M$  the two $\mathfrak{g}(\alpha)$-module structures on $M$ given by the maps $\xi$ and $d\mathfrak{J}$ differ by a character vanishing on the scalars. This character is given by the formula 
\begin{equation}\label{differencecharacter}
\xi-d\mathfrak{J}=\frac{1}{2}\sum_{i=0}^{\ell-1}(\alpha_{i-1}-\alpha_{i+1})\tr_{\mathfrak{gl}_{\alpha_i}}
\end{equation} 
where $\alpha_{-1}=0$.  In particular we have that, for every character $\chi$  of  $\mathfrak{g}(\alpha)$,  the $\chi$~-~equivariant subspace $M^{\mathfrak{g}(\alpha)}_{\chi}$ with respect to the $\mathfrak{g}(\alpha)$-action induced by the map $\xi$ is the same as the\  $\tilde{\chi}$~-~equivariant subspace with respect to the $\mathfrak{g}(\alpha)$-action  induced  by the map $d\mathfrak{J}$, where $\tilde{\chi}$ is the character defined by the following formula 
\begin{equation}\label{chitilde}
\tilde{\chi}:=\chi-\frac{1}{2}\sum_{i=0}^{\ell-1}(\alpha_{i-1}-\alpha_{i+1})\tr_{\mathfrak{gl}_{\alpha_i}}.   
\end{equation}
\end{prop}
\pf\ The proof is immediate. The fact that the  character in \eqref{differencecharacter}  vanishes on the scalars follows from the fact that both $\xi$ and $d\mathfrak{J}$  vanish on the scalars.
\epf   

Consider now the $W_{\alpha}$-module $M=S(\Rep_{\alpha}(Q^{op}))$, where the $W_{\alpha}$-action is as described in Remark \ref{Waction}. Let $\chi$ be a character of $\mathfrak{g}(\alpha)'$, and suppose that for $i=1,\dots, \ell-1$  there exist positive integers $b_i$  such that 
\begin{equation}\label{goodcharacter}
\chi_i=b_i+\frac{1}{2}(\alpha_{i-1}-\alpha_{i+1}).
\end{equation}
In this case we have
$$
\tilde{\chi}_i=\chi_i-\frac{1}{2}(\alpha_{i-1}-\alpha_{i+1})=b_i
$$
so that in particular the character $\tilde{\chi}$ is integral. By Theorem \ref{hamiltonianred} and Proposition \ref{mapsdifference} the $\chi$~-~equivariant subspace $M^{\mathfrak{g}(\alpha)'}_{\chi}$ of $M$ can be identified  with the space $\Gamma(\mathcal{L}_{\tilde{\chi}})$ of global sections on the flag variety $X/\Gg(\alpha)'\simeq\mathcal{F}\ell( \alpha_1,\dots, \alpha_{\ell-1}; N)$ of the line bundle 
$$
\mathcal{L}_{\tilde{\chi}}=X\times \C_{-\tilde{\chi}}/\Gg(\alpha)'
$$
where $\C_{-\tilde{\chi}}$ is the one-dimensional $\Gg(\alpha)'$-module corresponding to the character 
\begin{align*}
e^{-\tilde{\chi}}& =\prod_{i=1}^{\ell-1}\left(\dete_{\GL_{\alpha_i}}\right)^{-b_i}
\end{align*}
and the $\Gg(\alpha)'$-action on $X\times\C_{-\tilde{\chi}}$ is the diagonal one. Similarly to the case of the complete flag variety we have that if  $M^{\mathfrak{g}(\alpha)'}_{\chi}$ is nonzero then it is  a finite dimensional irreducible  $\GL_N$-module. We have the following proposition.
\begin{prop}\label{hstweight}
The space   $M^{\mathfrak{g}(\alpha)'}_{\chi}$ with the $\gln$-action induced by the map $\xi$  is an  irreducible  finite dimensional module  of highest weight 
\begin{equation}\label{hstformula}
\mu=\sum_{i=1}^{\ell-1}b_{i}\omega_{N-\alpha_i}-\left(\sum_{i=1}^{\ell-1} b_i+\frac{1}{2}\alpha_1\right)\omega_N.
\end{equation}
\end{prop}
\pf\ For $i=1,\dots, \ell-1$,  $r=1,\dots \alpha_i$, and $s=1,\dots,\alpha_{i-1}$ let us denote by $U^{\alpha_i,\alpha_{i-1}}_{r,s}$ the elementary $\alpha_i\times \alpha_{i-1}$ matrix whose $rs$-th entry is one, and whose remaining entries are zero.  For every choice of index $i$, these matrices clearly form a basis of $\Hom_{\C}(\C^{\alpha_{i-1}},\C^{\alpha_i})$. With this choice of bases, and using the same notation as in \eqref{goodbases}, we will identify  $\rho^{h^{\ast}_i, rs}$ with the representation of $Q^{op}$ of dimension vector $\alpha$ that assigns  the linear map $U^{\alpha_i,\alpha_{i-1}}_{r,s}$ to the arrow $h^{\ast}_i$, and assigns the zero map to all the other arrows.  These representations clearly form a basis of  $\Rep_{\alpha}(Q^{op})$, and we will call them the \emph{elementary representations} of $Q^{op}$ of dimension vector $\alpha$.

Let $\overline{\rho}^{h^{\ast}_i}$ be the  $\alpha_i\times \alpha_{i-1}$ matrix with entries in $\Rep_{\alpha}(Q^{op})\subset M$ given by the formula
$$
\overline{\rho}^{h^{\ast}_i}:=\left[\rho^{h^{\ast}_i, rs}\right]_{\alpha_i\times\alpha_{i-1}}
$$ 
In other words $\overline{\rho}^{h^{\ast}_i}$ is the matrix whose $rs$-th entry is the corresponding elementary representation of  $Q^{op}$ of dimension vector $\alpha$. 
 
 For any element $g=(g_0, g_1,\dots,g_{\ell-1})$  of $\GL_N\times \Gg(\alpha)'$, where $g_i\in \GL_{\alpha_i}$, we can then write the action of $g$ on $\Rep_{\alpha}(Q^{op})$ in matrix form as follows
$$
g\ast\overline{\rho}^{h^{\ast}_i}=\left[g_i\rho^{h^{\ast}_i, rs}g_{i-1}^{-1}\right]_{\alpha_i\times\alpha_{i-1}}
$$
where, as usual,  $g_i\rho^{h^{\ast}_i, rs}g_{i-1}^{-1}$ denotes the action of $\Gg(\alpha)$ on $\Rep_{\alpha}(Q^{op})$ by conjugation.

 If $g_i\in \GL(\alpha_i)$ and $g_{i-1}^{-1}\in \GL(\alpha_{i-1})$ are represented, in the standard bases of $\C^{\alpha_i}$ and $\C^{\alpha_{i-1}}$,  by the matrices $g_i=\left[g_{i,hk}^{\phantom{'}}\right]_{\alpha_i\times\alpha_i}$  and $g_{i-1}^{-1}=\left[g_{i-1, hk}'\right]_{\alpha_{i-1}\times\alpha_{i-1}}$ it is easy to calculate that
 $$
 g_i\rho^{h^{\ast}_i, rs}=\sum_{h=1}^{\alpha_i}g_{i,hr}\rho^{h^{\ast}_i, hs}\quad\mbox{ and }\quad \rho^{h^{\ast}_i, rs}g_{i-1}^{-1}=\sum_{h=1}^{\alpha_{i-1}}\rho^{h^{\ast}_i, rh}g_{i-1,sh}'.
 $$
 From these  formulas we can easily deduce  the following convenient description of the action of $\GL_N\times\Gg(\alpha)'$ on $\Rep_{\alpha}(Q^{op})$ in matrix form 
\begin{equation}\label{matrixaction}
g\ast\overline{\rho}^{h^{\ast}_i}=g_i^T\overline{\rho}^{h^{\ast}_i}(g_{i-1}^T)^{-1}
\end{equation}
where $g_i^T$ and $(g_{i-1}^T)^{-1}$ are the transpose matrices of $g_i$ and $g_{i-1}^{-1}$ respectively,  and $g_i^T\overline{\rho}^{h^{\ast}_i}(g_{i-1}^T)^{-1}$ on the right hand side of the
equality is  just the ordinary product of the corresponding matrices.

For every $i=1,\dots,\ell-1$ we can regard the product
$$
\pi_i=\overline{\rho}^{h^{\ast}_i}\dots \overline{\rho}^{h^{\ast}_1}
$$
as an $\alpha_i\times N$-matrix with entries in $S(\Rep_{\alpha}(Q^{op}))=M$.
Let us denote by   $D_i$ the leading principal minor of size $\alpha_i$  of the matrix $\pi_i$, and let us define the vector $v$ in $M$  as follows
$$
v:=D_1^{b_1}\dots D_{\ell-1}^{b_{\ell-1}}.
$$ 
We will show  that $v$ is in $M^{\mathfrak{g}(\alpha)'}_{\chi}$, and that it is a  lowest weight  vector  for the $\gln$-action induced by the map $d\mathfrak{J}$.

First of all, observe  that, using formula \eqref{matrixaction},  for any $g'=(g_1,\dots,g_{\ell-1})$ in $\Gg(\alpha)'$  and  any $i=1,\dots,\ell-1$ we have
\begin{align*}
g'\ast\pi_i&=(g'\ast\overline{\rho}^{h^{\ast}_i})(g'\ast\overline{\rho}^{h^{\ast}_{i-1}})\cdots(g'\ast\overline{\rho}^{h^{\ast}_1})\\
&=g_i^T\overline{\rho}^{h^{\ast}_i}(g_{i-1}^T)^{-1}g_{i-1}^T\overline{\rho}^{h^{\ast}_{i-1}}(g_{i-2}^T)^{-1}\cdots g_{\ell-1}^T\overline{\rho}^{h^{\ast}_1}\\
            &=g_i^T\pi_i.
\end{align*}
Now since $\pi_i$ is an $\alpha_i\times N$ matrix we can write it in block form as follows
$$
\pi_i=\left(A_i|B_i\right)
$$
where $A_i$ is of size $\alpha_i\times\alpha_i$ and $B_i$ is of size $\alpha_i\times(N-\alpha_i)$.  With this notation  we have that $D_i=\det(A_i)$ and  
$$
g'\ast \pi_i=g_i^T\pi_i=\left(g_i^TA_i|g_i^TB_i\right).
$$
Thus, calculating the action of $g'$ on $D_i^{b_i}$, we get
\begin{align*}
g'\ast D_i^{b_i}&=(g'\ast D_i)^{b_i}\\&=(\det(g_i^TA_i))^{b_i}\\&=(\det g_i^T)^{b_i}\det(A_i)^{b_i}\\&=(\det g_i)^{b_i}D^{b_i}_i.
\end{align*}
It follows that
\begin{align*}
g'\ast v&=\prod_{i=1}^{\ell-1}(g'\ast D_i)^{b_i}\\
        &=\prod_{i=1}^{\ell-1}(\det g_i)^{b_i}D_i^{b_i}\\
        &=\left(\prod_{i=1}^{\ell-1}(\det g_i)^{b_i} \right)v.
\end{align*}
Thus $v$ is in $\Gamma(\mathcal{L}_{\tilde{\chi}})=M^{\mathfrak{g}(\alpha)'}_{\chi}$, and $M^{\mathfrak{g}(\alpha)'}_{\chi}$ is nonzero.

Let now  $T$ be the maximal torus of diagonal matrices in $\mathrm{GL}_N$, and let $t_0=\mathrm{diag}(c_1,\dots,c_N)$ be in $T$. We have
\begin{align*}
t_0\ast \pi_i&=\pi_i(t_0^T)^{-1}\\
              &=\pi_it_0^{-1}\\
              &=(A_i|B_i)t_0^{-1}\\
              &=\left(A_i\mathrm{diag}(c_1^{-1},\dots,c_{\alpha_i}^{-1})|B_i\mathrm{diag}(c_{\alpha_i+1}^{-1},\dots,c_N^{-1})\right).
\end{align*}
Calculating the action of $t_0$ on $D_i^{b_i}$, we get
\begin{align*}
t_0\ast D_i^{b_i}&=(t_0\ast D_i)^{b_i}\\
                  &=\det(A_i\mathrm{diag}(c_1^{-1},\dots,c_{\alpha_i}^{-1}))^{b_i}\\
                  &=\left(c_1^{-1}\cdots c_{\alpha_i}^{-1}\right)^{b_i}\det(A_i)^{b_i}\\
                  &=\left(c_1\cdots c_{\alpha_i}\right)^{-b_i}D_i^{b_i}.
\end{align*}
This gives  the following formula for the action of $t_0$ on $v$
$$
t_0\ast v=\left(\prod_{i=1}^{\ell-1}\left(c_1\cdots c_{\alpha_i}\right)^{-b_i}\right) v =e^{\mu'}(t_0)v
$$
where $\mu'$ is the following weight
$$
\mu'=-\sum_{i=1}^{\ell-1}b_i\omega_{\alpha_i}.
$$
Therefore $v$ is a  weight vector of weight $\mu'$ for the action of $\gln$ induced by $d\mathfrak{J}$.

Finally, let $B$ be the Borel subgroup of upper triangular matrices in $\mathrm{GL}_N$, and let $B'$ be the opposite Borel of lower triangular matrices. Let $U'$ be the unipotent radical of $B'$, and let $l_0$ be an element of $U'$, that is $l_0$ is a unipotent lower triangular matrix in $\GL_N$.
 For every $i=1,\dots,\ell-1$ we have
\begin{align*}
l_0\ast\pi_i&=\pi_i(l_0^T)^{-1}                    
\end{align*}
Here the transpose matrix  $(l_0^T)^{-1}$ is a unipotent \emph{upper} triangular matrix, and can  thus be written in block form as 
$$
(l_0^T)^{-1}=\left(\begin{array}{cc} U_{\alpha_i}& \ast\\ 0 &  U_{N-\alpha_i}\end{array}\right)
$$
where $U_{\alpha_i}$ and $U_{N-\alpha_i}$ are unipotent upper triangular matrices  in $\mathrm{GL}_{\alpha_i}$ and $\mathrm{GL}_{N-\alpha_i}$ respectively,  and $0$ stands for the zero matrix of size $(N-\alpha_i)\times \alpha_i$.
 Thus we have
 \begin{align*}
l_0\ast\pi_i&=\pi_i(l_0^T)^{-1} \\
            &=(A_i|B_i)\left(\begin{array}{cc}U_{\alpha_i}& \ast\\ 0 & U_{N-\alpha_i}\end{array}\right) \\
           &=\left(A_iU_{\alpha_i}|C+B_iU_{N-\alpha_i}\right)                 
\end{align*}
where $C$ is some $\alpha_i\times(N-\alpha_i)$ matrix. Calculating the action of $l_0$ on $D_i^{b_i}$ we get
\begin{align*}
l_0\ast D_i^{b_i}&=(l_0\ast D_i)^{b_i}\\&=(\det(A_iU_{\alpha_i}))^{b_i}\\&=\det(A_i)^{b_i}\det(U_{\alpha_i})^{b_i}\\&=\det(A_i)^{b_i}\\&=D_i^{b_i}.
\end{align*}
where we used the fact that $\det(U_{\alpha_i})$ is equal to one. Therefore 
$$
l_0\ast v=v
$$
and $v$ is a lowest weight vector for the corresponding $\gln$ action.
 
Thus $M^{\mathfrak{g}(\alpha)'}_{\chi}$ is an irreducible module of lowest weight $\mu'$ for  the $\gln$-action induced by the map $d\mathfrak{J}$. It now follows from Proposition \ref{mapsdifference}  that, with respect to the $\gln$-action induced by the map $\xi$, the space $M^{\mathfrak{g}(\alpha)'}_{\chi}$ is a finite dimensional irreducible module of lowest weight
$$
\mu''=\mu'-\frac{1}{2}\alpha_1\tr_{\gln}=-\sum_{i=1}^{\ell-1}b_i\omega_{\alpha_i}-\frac{1}{2}\alpha_1\omega_N.
$$

To obtain the highest weight from the lowest weight we just have to apply the transformation $-\tau$ to $\mu''$, where  $\tau$ is the canonical involution for type $A$. Since   $-\tau$ is given by the formulas 
\begin{align*}
-\tau(\omega_N)&=\omega_N \\
-\tau(\omega_r)&=\omega_N-\omega_{N-r},\quad\quad r=1,\dots,N-1  
\end{align*}
we get the desired expression for $\mu$.
\epf

\section{Comparison of representations}\label{comparison}
 With the results of the previous section at hand, we are now ready to show how to obtain some of the representations of Section \ref{gdahasrep} as a special case of the representations of Section \ref{ggrep}, using the isomorphism $\varphi$ of Proposition \ref{isomorphism}.

Let  $D$ be a star-shaped graph that is not a  finite Dynkin diagram. Let the vertices of $D$ be labeled as in Section \ref{gdahadefisec}. Let $Q$ be the quiver obtained from $D$ by assigning to the edges along each leg the orientation toward the node.
Let $N$ be a positive integer. Let $\alpha$ be a dimension vector for $Q$ such that $\alpha_{0}=N$, and the dimensions $\alpha_{(k,s)}$ at the remaining vertices are all nonzero and   strictly decreasing  along each leg, as one moves away from the node. In other words for $k=1,\dots, m$ we have
\begin{equation}\label{flagsequence2}
\alpha_{0}=N>\alpha_{(k,1)}>\cdots>\alpha_{(k,\ell_k-1)}>0 .
\end{equation}  
Denote by $Q_k$ the quiver (of type $A$) corresponding to the $k$-th leg of  $Q$ (including the node), and let $\alpha^{(k)}$ be the restriction of the dimension vector $\alpha$ to $Q_k$. Let $E_k~:=~\Rep_{\alpha^{(k)}}(\overline{Q}_k)~\subset~\Rep_{\alpha}(\overline{Q})$ be the corresponding symplectic subspace, and denote by $W_k\subset W_{\alpha}$ its Weyl algebra. Then it is clear that  the vector space $\Rep_{\alpha}(\overline{Q})$ has the following decomposition
$$
\Rep_{\alpha}(\overline{Q})=\bigoplus_{k=1}^m E_{k}.
$$
It follows that the Weyl algebra $W_{\alpha}$ is the tensor product
\begin{equation}\label{decoweyl}
W_{\alpha}=\bigotimes_{k=1}^m W_{k}.
\end{equation}
To ease the notation let us set 
$$
\Gg_{k}:=\Gg(\alpha^{(k)})'=\prod_{s=1}^{\ell_k-1}\GL_{\alpha_{(k,s)}}\qquad \mbox{ and }\qquad  \mathfrak{g}_k:=\mathfrak{g}(\alpha^{(k)})'=\prod_{s=1}^{\ell_k-1}\mathfrak{gl}_{\alpha_{(k,s)}}.
$$ 
We  have
\begin{equation}
\Gg(\alpha)=\GL_N\times\prod_{k=1}^m \Gg_{k},\qquad \qquad \mathfrak{g}(\alpha)=\gln\times\prod_{k=1}^m \mathfrak{g}_{k}.
\end{equation}  

Let us now consider the  construction of Section \ref{ggrep} in this special case. For any $W_{\alpha}$-module $M$  we can write  $M=\bigotimes_{k=1}^mM_k$, where $M_k$ is a $W_k$-module. Let  $\chi$ be a character of $\mathfrak{g}(\alpha)$. Let us write 
\begin{align*}
\chi'&=\chi_{0} \tr_{\mathfrak{gl}_N}=c \tr_{\mathfrak{gl}_N}\\
\chik &=\sum_{s=1}^{\ell_k-1}\chi_{(k,s)} \tr_{\mathfrak{gl}_{\alpha_{(k,s)}}}
\end{align*}
for the restriction of $\chi$ to $\mathfrak{gl}_N$ and to $\mathfrak{g}_k$ respectively.  We can now construct the module $F_{n,\chi}(M)$ in two steps, first by taking the $\sum_{k=1}^{m}\chik$-equivariant space with respect to the $\prod_{k=1}^{m}\mathfrak{g}_k$-action, and then taking the $\chi'$-equivariant space with respect to the  $\gln$-action. By restriction, we get a representation of  the spherical subalgebra $e_{0}^{\otimes n}\mathcal{A}_{n,\lambda,\nu}(Q)e_{0}^{\otimes n}$ on the space 
\begin{align*}
e_{0}^{\otimes n}F_{n,\chi}(M)&=e_{0}^{\otimes n}(M\otimes U^{\otimes n})^{\mathfrak{g}(\alpha)}_{\chi}\\
&=(M\otimes e_{0}^{\otimes n}U^{\otimes n})^{\mathfrak{g}(\alpha)}_{\chi}\\
&=\left(M\otimes (\C^N)^{\otimes n}\right)^{(\gln\times\prod_{k=1}^m\mathfrak{g}_k)}_{\chi}\\
&=\left(\bigotimes_{k=1}^m(M_k)^{\mathfrak{g}_k}_{\chik}\otimes (\C^N)^{\otimes n}\right)_{\chi'}^{\gln}.
\end{align*}

Using  the results of Section \ref{borelweil}, we will show how, for appropriate choices of the $W_{\alpha}$-module $M$ and of the character $\chi$,  the subspace $(M)_{\chik}^{\tilde{\mathfrak{g}}_k}$ is a finite dimensional $\ell_k$-stepped $\gln$-module for every $k=1,\dots,m$, and the above representation is of the type  described in Theorem \ref{repreBn}.

Suppose  there exist positive integers $b_{ks}$  such that for $k=1,\dots,m$, $s=1,\dots,\ell_k-1$ we have
\begin{equation}\label{characterspecial}
\chi_{(k,s)}=b_{ks}+\frac{1}{2}(\alpha_{(k,s-1)}-\alpha_{(k,s+1)}).
\end{equation}
Suppose moreover that $\chi$ satisfies condition \eqref{conditionchis}, that is we have the identity
$$
\sum_{k=1}^m\sum_{s=1}^{\ell_k-1}\chi_{(k,s)}\alpha_{(k,s)}+cN-n=0.
$$
Substituting the values for $\chi_{(k,s)}$ given in formula \eqref{characterspecial} in the above identity we get the following condition on the integers $b_{k,s}$ and the complex number $c$ that define the character $\chi$
\begin{equation}\label{conditionchisstar}
\sum_{k=1}^m\sum_{s=1}^{\ell_k-1}b_{k,s}\alpha_{k,s}+\frac{1}{2}\sum_{k=1}^m\alpha_{(k,1)}N+cN-n=0.
\end{equation}

Set $M_k=S(\Rep_{\alpha^{(k)}}(Q^{op}_k))$, where the $W_{k}$-action is as described in Remark \ref{Waction}. 

As in Section \ref{flagsection} define the open subset $X_k=\oplus_{s=1}^{\ell-1}\mathrm{Incl}(\C^{\alpha_{(k,s)}}, \C^{\alpha_{(k,s-1)}})$ of $\Rep_{\alpha_k}(Q_k)$. Then, by Theorem \ref{hamiltonianred} and Proposition \ref{mapsdifference} applied to  the algebra $W_k$ and the character $\chik$, the  subspace $(M_k)^{\mathfrak{g}_k}_{\chik}$ can be identified  with the space $\Gamma(\mathcal{L}_{\chikti})$ of global sections on the flag variety $X_k/\Gg_k\simeq\mathcal{F}\ell( \alpha_{(k,1)},\dots, \alpha_{(k,\ell_k-1)}; N)$ of the line bundle 
$$
\mathcal{L}_{\chikti}=X_k\times \C_{-\chikti}/\Gg_k,
$$
where $\C_{-\chikti}$ is the one-dimensional $\Gg_k$-module corresponding to the character 
\begin{align*}
-\chikti&=\prod_{s=1}^{\ell_k-1}\left(\dete_{\GL_{\alpha_{(k,s)}}}\right)^{-(\chi_{(k,s)}-\frac{1}{2}(\alpha_{(k,s-1)}-\alpha_{(k,s+1)}))}\\
             &=\prod_{s=1}^{\ell_k-1}\left(\dete_{\GL_{\alpha_{(k,s)}}}\right)^{-b_{ks}}
\end{align*}
and the $\Gg_k$-action on $X_k\times\C_{-\chikti}$ is the diagonal one. By Proposition \ref{hstweight}, the space  $(M_k)^{\mathfrak{g}_k}_{\chik}$   is an  irreducible  finite dimensional $\gln$-module  of highest weight 
\begin{equation}\label{highweightformula}
\mu_k=\sum_{s=1}^{\ell_k-1}b_{ks}\omega_{N-\alpha_{(k,s)}}-\left(\sum_{s=1}^{\ell_{k}-1} b_{ks}+\frac{1}{2}\alpha_{(k,1)}\right)\omega_N.
\end{equation}

It follows that choosing $\displaystyle M=S(\Rep_{\alpha}(Q^{op}))=\otimes_{k=1}^{m}S(\Rep_{\alpha^{(k)}}(Q^{op}_k))=\otimes_{k=1}^mM_k$ we get a representation of the algebra 
$e_{0}^{\otimes n}\mathcal{A}_{n, \lambda, -1}(Q)e_{0}^{\otimes n}$ on the space
\begin{align*}
e_{0}^{\otimes n}F_{n,\chi}(M)&=\left(\bigotimes_{k=1}^mV_k\otimes (\C^N)^{\otimes n}\right)^{\gln}_{\chi'}=E_{n,\chi'}(V_1,\dots, V_m)
\end{align*}
where $V_k$ is the irreducible $\gln$-module whose highest weight $\mu_k$ is described by  formula (\ref{highweightformula}). 
We would like to show that this representation is of the type described in Theorem \ref{repreBn}. 

First of all, observe that, with our choice of character $\chi$ (see (\ref{characterspecial})) and for $\nu=-1$, according to Theorem \ref{representationtheorem}, the parameter $\lambda$ must satisfy 
\begin{equation}\label{lambda02}
\lambda_{0}=c-N+\frac{1}{2}\sum_{k=1}^m\alpha_{(k,1)}
\end{equation}
\begin{equation}\label{lambdas2}
\lambda_{(k,s)}=b_{ks}-\alpha_{(k,s)}+\alpha_{(k,s-1)}
\end{equation}
Note that, since the numbers $\alpha_{(k,s)}$ decrease along the $k$-th leg moving away from the node,  the numbers $N-\alpha_{(k,s)}$  increase. With  the same notation as in Section \ref{gdahasrep}, we can thus write 
$$
\mu_k=b_{k1}\omega_{p_{k1}}+\dots+b_{k,\ell_k-1}\omega_{p_{k(\ell_{k}-1)}}+b_{k\ell_k}\omega_N,
$$
where  $p_{ks}=N-\alpha_{(k,s)}$ and   $b_{ks}$ are the positive integers in 
(\ref{characterspecial}) for $s=1,\dots, \ell_k-1$, while  $b_{k\ell_k}=-\left(\sum_{s=1}^{\ell_{k}-1} b_{ks}+\frac{1}{2}\alpha_{(k,1)}\right)$.  It is now easy to calculate that  the value of the parameter  $\lambda$ given in formulas (\ref{lambda02}) and (\ref{lambdas2}) agrees with the value given in formulas (\ref{lambda0}), (\ref{lambdas}) in Section \ref{gdahasrep} for $\nu=1$.
Thus  choosing $\xi_k$ as in formula (\ref{xik}), that is 
\begin{equation}
\xi_k=1\cdot\left(\frac{1}{m}\sum_{h\neq k}\sum_{j=1}^{\ell_h}b_{hj}+\frac{1-m}{m}\sum_{j=1}^{\ell_k}b_{ks}\right)=-\frac{1}{2m}\sum_{h\neq k}\alpha_{(h,1)}-\frac{1-m}{2m}\alpha_{(k,1)}
\end{equation}
we have an isomorphism  $\varphi: B_n(\g,1)\simeq e_{0}^{\otimes n}\mathcal{A}_{n, \lambda, -1}(Q)e_{0}^{\otimes n}$ as in Theorem \ref{isomorphism}, where $\g$ is as in Lemma \ref{Bneigenvalues}.

Moreover we can check that the condition  on the character $\chi$ in formula \eqref{conditionchisstar} agrees with the condition on the character $\chi'=c\tr_{\mathfrak{gl}_N}$  in formula \eqref{traceconditionBn}  when $\mu_k$ is as in  \eqref{highweightformula}. Indeed, since for every $k$ we have $b_{k\ell_k}=-\left(\sum_{s=1}^{\ell_{k}-1} b_{ks}+\frac{1}{2}\alpha_{(k,1)}\right)$,  according to formula \eqref{traceconditionBn} we must have
$$
\sum_{k=1}^m\sum_{s=1}^{\ell_k-1}b_{ks}(N-\alpha_{(k,s)})-\sum_{k=1}^m\left(\sum_{s=1}^{\ell_{k}-1} b_{ks}+\frac{1}{2}\alpha_{(k,1)}\right)N+n-cN=0
$$
and so
$$
\sum_{k=1}^m\sum_{s=1}^{\ell_k-1}b_{ks}\alpha_{(k,s)}+\frac{1}{2}\sum_{k=1}^m\alpha_{(k,1)}N+cN-n=0,
$$
which is exactly the identity in \eqref{conditionchisstar}.

We now want to prove that pulling back the above representation  of $e_{0}^{\otimes n}\mathcal{A}_{n, \lambda, -1}(Q)e_{0}^{\otimes n}$ on the space  $E_{n,\chi'}(V_1,\dots,V_m)$ via the isomorphism $\varphi$ we get the representation of the algebra $B_n(1,\g)$ in Theorem \ref{repreBn}.

To ease the notation  set $W=E_{n,\chi'}(V_1,\dots,V_m)$ and let $\psi:e_{0}^{\otimes n}\mathcal{A}_{n, \lambda, -1}(Q)e_{0}^{\otimes n}\rightarrow \mathrm{End}(W)$ be the homomorphism corresponding to the above representation. Let $\U_{rs}$ denote the $rs$-th matrix unit in $\mathfrak{gl}_{\alpha_0}=\gln$. Using once more the pair of dual bases for $\Rep_{\alpha}(\overline{Q})$ defined by formulas (\ref{goodbases}) we have
\begin{align*}
\psi(\varphi(Y_{i,k}))&=\psi\left(e_{0}^{\otimes i-1}\otimes\left(-h_{k1}h_{k1}^{\ast}+\left(\xi_k+\frac{\lambda_{0}}{m}\right)e_{0}\right)\otimes e_{0}^{\otimes n-i}\right)\\
&=-\sum_{p,q}\rho_p\rho_q|_{M}\rho^p(h_{k1})\rho^q(h^{\ast}_{k1})|_{i+m}-\frac{(N-c)}{m}+\frac{1}{2}\alpha_{(k,1)}\\
&=-\!\!\!\!\!\!\!\!\!\!\!\!\sum_{\begin{array}{c}\scriptscriptstyle{s,r =1,\dots,N}\\\scriptscriptstyle{t,l=1,\dots,\alpha_{(k,1)}}\end{array}}\!\!\!\!\!\!\!\!\!\left(\rho_{h^{\ast}_{k1},ts}\rho_{h_{k1},rl}\right)|_{M_k}\left(\rho^{h_{k1},st}(h_{k1})\rho^{h^{\ast}_{k1},lr}(h^{\ast}_{k1})\right)|_{i+m}-\frac{N-c}{m}+\frac{1}{2}\alpha_{(k,1)}\\
&=-\!\!\!\!\!\!\!\!\!\!\!\!\sum_{\begin{array}{c}\scriptscriptstyle{s,r =1,\dots,N}\\\scriptscriptstyle{t,l=1,\dots,\alpha_{(k,1)}}\end{array}}\!\!\!\!\!\!\!\!\!\delta_{t,l}\left(\rho_{h^{\ast}_{k1},ts}\rho_{h_{k1},rl}\right)|_{M_k}\left(-\U_{sr}\right)|_{i+m}-\frac{N-c}{m}+\frac{1}{2}\alpha_{(k,1)}\\
&=\!\!\!\!\!\!\sum_{\begin{array}{c}\scriptscriptstyle{s,r =1,\dots,N}\\\scriptscriptstyle{t=1,\dots,\alpha_{(k,1)}}\end{array}}\!\!\!\!\!\!\!\!\!\left(\rho^{h^{\ast}_{k1},ts}\rho_{h_{k1},rt}\right)|_{M_k}\,\U_{sr}|_{i+m}-\frac{N-c}{m}+\frac{1}{2}\alpha_{(k,1)}\\
&=\!\!\!\!\!\sum_{\begin{array}{c}\scriptscriptstyle{s,r =1,\dots,N}\end{array}}\!\!\!\!\!\!\!\!\!\left(-\xi(\U_{rs})-\delta_{rs}\frac{1}{2}\alpha_{(k,1)}\right)|_{M_k}\,\U_{sr}|_{i+m}-\frac{N-c}{m}+\frac{1}{2}\alpha_{(k,1)}\\
&=-\!\!\!\!\!\!\!\!\!\!\sum_{\begin{array}{c}\scriptscriptstyle{s,r =1,\dots,N}\end{array}}\!\!\!\!\!\!\!\!\!
\xi(\U_{rs})|_{M_k}\,\U_{sr}|_{i+m}-\frac{N-c}{m}\\
&=-\Omega_{k,i+m}-\frac{N-c}{m}
\end{align*}
where we used  the fact $\rho_{h^{\ast}_{k1},ts}=\rho^{h^{\ast}_{k1},ts}$, and we used formula (\ref{xiformula}) to compute $\xi(\U_{rs})$. Note that we get exactly the same formula as in Theorem \ref{repreBn} for $\nu=1$.

\end{document}